\documentclass[11pt, reqno]{amsart}
\usepackage{amsthm,amsmath,amssymb}
\usepackage{mathtools}
\usepackage{graphicx}
\usepackage{caption}
\usepackage{ragged2e}
\graphicspath{{Figure/}}
\numberwithin{equation}{section}
\numberwithin{figure}{section}
\theoremstyle{plain}

\newtheorem{thm}{Theorem}[section]

\newtheorem{cor}[thm]{Corollary}

\theoremstyle{definition}
\theoremstyle{remark}

\begin{document}
\setlength{\abovedisplayskip}{10pt}
\setlength{\belowdisplayskip}{10pt}

\title{Lozenge tilings of a hexagon with a horizontal intrusion}

\author{Seok Hyun Byun}
\thanks{This research was supported in part by Lilly Endowment, Inc., through its support for the Indiana University Pervasive Technology Institute.}

\address{Department of Mathematics, Indiana University, Bloomington}
\email{byunse@indiana.edu}

\maketitle

\begin{abstract}
Motivated by a conjecture posed by Fulmek and Krattenthaler, we provide product formulas for the number of lozenge tilings of a semiregular hexagon containing a horizontal intrusion. As a direct corollary, we obtain a product formula for the number of boxed plane partitions with a certain restriction. We also investigate the asymptotic behavior of the ratio between the number of lozenge tilings of a semiregular hexagon containing a horizontal intrusion and that of a semiregular hexagon without an intrusion.
\end{abstract}

\section{Introduction}
On the unit triangular lattice where one of the lattice lines is vertical, we draw a hexagon whose parallel sides have the same side lengths. Such a hexagon is \textit{a semiregular hexagon}. One natural question to ask is how many lozenge tilings\footnote{A \textit{lozenge} is a union of two adjacent unit triangles. A \textit{lozenge tiling} of a region is a collection of lozenges that cover the region without any gaps or overlaps. For example, see the left picture in Figure 2.2 or the left picture in Figure 2.4.} does this region have. In fact, the answer to this question was given a long time ago. The bijection of David and Tomei [9] allows one to interpret MacMahon's classical theorem for boxed plane partitions [25, \S 429, $q\to 1$; proof in \S 494] as follows: the number of lozenge tilings of a hexagon with sides of length \(a\), \(b\), \(c\), \(a\), \(b\), \(c\) (clockwise from the left) is given by 
\begin{equation}
    \prod_{i=1}^{a}\prod_{j=1}^{b}\prod_{k=1}^{c}\frac{i+j+k-1}{i+j+k-2}=\frac{H(a)H(b)H(c)H(a+b+c)}{H(a+b)H(b+c)H(c+a)},
\end{equation}
where $H(n):=\displaystyle \prod_{i=0}^{n-1} i!$.

To answer and to generalize the problems posed by Propp (Problem 1 and Problems 4 in [27]), many authors considered the enumeration of lozenge tilings of a hexagon with a fixed lozenge (for examples, see [5], [10], [12], [13], and [14]). The motivation for this paper is a conjecture posed by Fulmek and Krattenthaler in [12] (it is Conjecture 1 of their paper). In [12], they considered a hexagon with sides of length $M, N, N, M, N, N$ (clockwise from the left) and enumerated the number of lozenge tilings with a fixed lozenge on the horizontal symmetry axis (see the left picture in Figure 1.1 for an example). In that paper, they conjectured various formulas that enumerate the number of lozenge tilings of the same hexagon containing several lozenges fixed on the axis. In particular, when the fixed lozenges are left-aligned on the axis, the conjectured formula (equation (4.4) in Conjecture 1 of their paper) was ``nice" in the sense that each factor in the formula is linear in each parameter. More precisely, for nonnegative integer $m$ and positive integers $N$ and $r$ with $N\geq r$, they considered a hexagon with sides of length $2m, N, N, 2m, N, N$ (clockwise from the left) and its lozenge tilings. They conjectured that the number of lozenge tilings of the hexagon which contains $r$ leftmost lozenges on the horizontal symmetry axis (see the middle picture in Figure 1.1 for an example) is given by\footnote{For a positive integer $n$, 
$n!!:=
\begin{cases}
    n\cdot(n-2)\cdots4\cdot2  & \text{if $n$ is even,}\\
    n\cdot(n-2)\cdots3\cdot1  & \text{if $n$ is odd.}
\end{cases}$
}
\begin{equation}
\begin{aligned}
    2^{\frac{1}{2}(r-1)(r-2N)}\frac{\binom{m+N-1}{m}^2}{\binom{2m+2N-1}{2m}}\prod_{i=N-r}^{N-2}\frac{1}{i!}\prod_{i=1}^{r-1}&\frac{(2i)!!(2N-2i-1)!!(m+i+1)_{N-2i-1}}{(2i-1)!!(m+i+\frac{1}{2})_{N-2i}}\\
    &\times\prod_{i=1}^{N}\prod_{j=1}^{N}\prod_{k=1}^{2m}\frac{i+j+k-1}{i+j+k-2},
\end{aligned}
\end{equation}
where the \textit{shifted factorial} $(a)_n$ is $(a)_n:=a(a+1)\cdots(a+n-1)$ for a positive integer $n$ and $(a)_0:=1$.

This conjecture was already proven by Ciucu and Krattenthaler. In [6], they provided a sketch of the proof of this conjecture. Their proof was based on the Matching Factorization Theorem of Ciucu in [1] and the result of Ciucu in [2] that enumerated the number of lozenge tilings of certain families of regions.
In this paper, we generalize the above result to arbitrary semiregular hexagons (see the right picture in Figure 1.1 for an example). More precisely, we enumerate the number of lozenge tilings of arbitrary semiregular hexagons with even numbers of left-aligned unit triangles on the perpendicular bisector of the left side removed. As a direct corollary, we obtain a product formula for the number of boxed plane partitions with a certain restriction. We also analyze the asymptotic behavior of the formulas in a scaling limit.

\begin{figure}
    \centering
    \includegraphics[width=16cm]{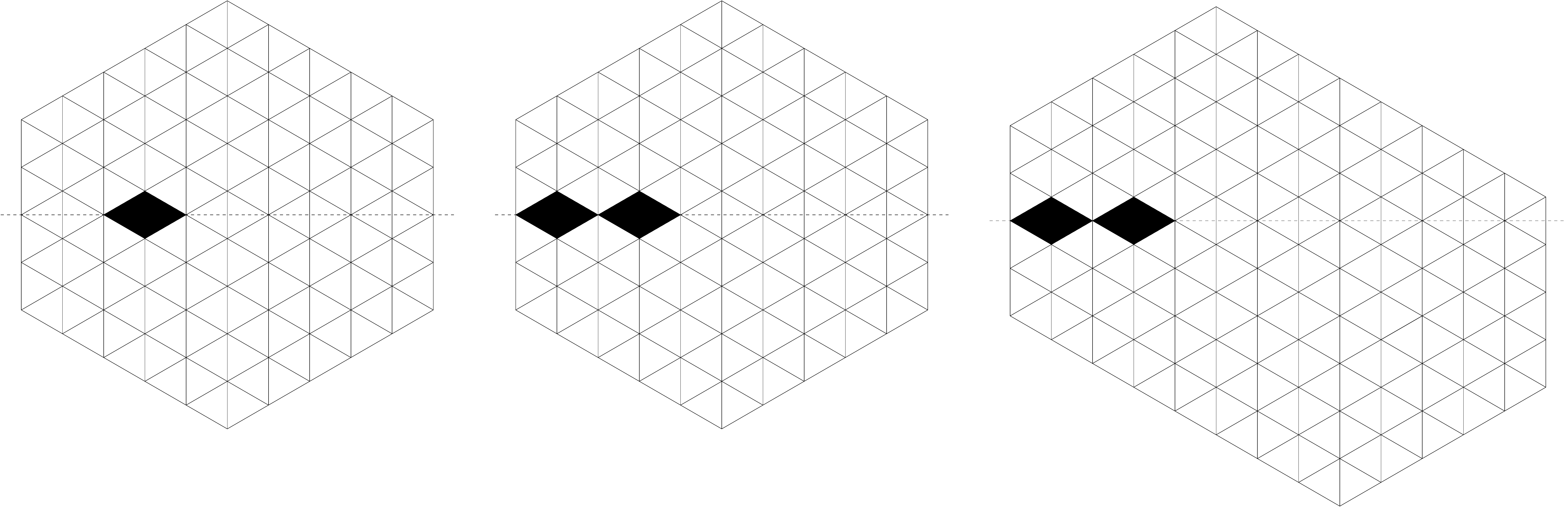}
    \caption{The left picture is one of the regions studied by Fulmek and Krattenthaler. The middle picture is one of the regions whose number of lozenge tilings was conjectured by them. The right picture is one of the regions studied in the current paper.}
\end{figure}

This paper is organized as follows:
in Section 2, we state our main results.
In Section 3, we provide several theorems from the literature that we use in Section 4.
In Section 4, we give a proof of the main theorem, Theorem 2.1.
Finally, in Section 5, we prove Theorem 2.3, which concerns the asymptotic behavior of the formulas obtained in Theorem 2.1.

\section{Statement of main results}

We consider a triangular lattice, comprised of unit equilateral triangles where one of the lattice lines is vertical. For nonnegative integers $a, b,$ and $c$, consider a hexagon with sides of length $a$, $b$, $c$, $a$, $b$, and $c$ (clockwise from the left) on the lattice and denote it by $H_{a,b,c}$.
For nonnegative integers $a, b, c,$ and $d$ such that $d\leq min(b,c)$, we define two regions $H_{2a,b,c;d}$ and $H_{2a+1,b,c;d}$ (see Figure 2.1 for two examples).

When $d=0$, we set $H_{2a,b,c;0}:=H_{2a,b,c}$ and $H_{2a+1,b,c;0}:=H_{2a+1,b,c}$. When $d>0$, the region $H_{2a,b,c;d}$ is obtained from the hexagonal region $H_{2a,b,c}$ by deleting $2d$ unit triangles as follows: consider the perpendicular bisector $l$ of the left side of the region $H_{2a,b,c}$. Label the unit triangles on $l$ that are contained in $H_{2a,b,c}$ by $1, 2, ...$ from the left. $H_{2a,b,c;d}$ is the region obtained from $H_{2a,b,c}$ by deleting $2d$ unit triangles labeled by $1, 2, ... , 2d$. One can see these deleted $2d$ unit triangles form $d$ consecutive unit lozenges. The region $H_{2a+1,b,c;d}$ is defined in a similar way. In the case of $H_{2a+1,b,c;d}$, one can see the deleted $2d$ unit triangles form $d$ consecutive unit bowties.
Two regions $H_{2a,b,c;d}$ and $H_{2a+1,b,c;d}$ are well-defined for any nonnegative integers $a$, $b$, $c$, and $d$ if all $2d$ unit triangular holes are contained in the hexagons $H_{2a,b,c}$ and $H_{2a+1,b,c}$, respectively. However, in this paper, we only deal with the case when the parameters satisfy $d \leq b \leq c$.
The reason why we assume the first inequality $d \leq b$ is that if $d > b$, then the region $H_{2a,b,c;d}$ (and similarly $H_{2a+1,b,c;d}$) has no lozenge tiling. One way to see this is using the nonintersecting lattice path interpretation of lozenge tilings. As one can see on the left picture in Figure 2.2, by selecting suitable start and end points on the hexagon, each tiling corresponds to a family of nonintersecting paths across lozenges (in that picture, start and end points are marked by red and blue dots). The paths can then be translated into a family of nonintersecting paths on the integer lattice. The number of tilings of the region is thus the number of families of nonintersecting paths between start and end points on the integer lattice where paths can only move along lattice lines. If $d>b$, as shown in the 
right picture in Figure 2.2, there are more end points on the holes than there are start points on the northwest side. This means there can be no family of nonintersecting paths that connect all start points to all end points, thus the region has no tilings in that case.
Also, one can easily see that two regions $H_{2a,b,c;d}$ and $H_{2a,c,b;d}$ (and similarly $H_{2a+1,b,c;d}$ and $H_{2a+1,c,b;d}$) have the same number of lozenge tilings since $H_{2a,c,b;d}$ can be obtained by reflecting $H_{2a,b,c;d}$ across the perpendicular bisector $l$. Hence, we can restrict our attention to the case when $b \leq c$, and the formulas for the case when $c \leq b$ follow.

\begin{figure}
    \centering
    \includegraphics[width=12cm]{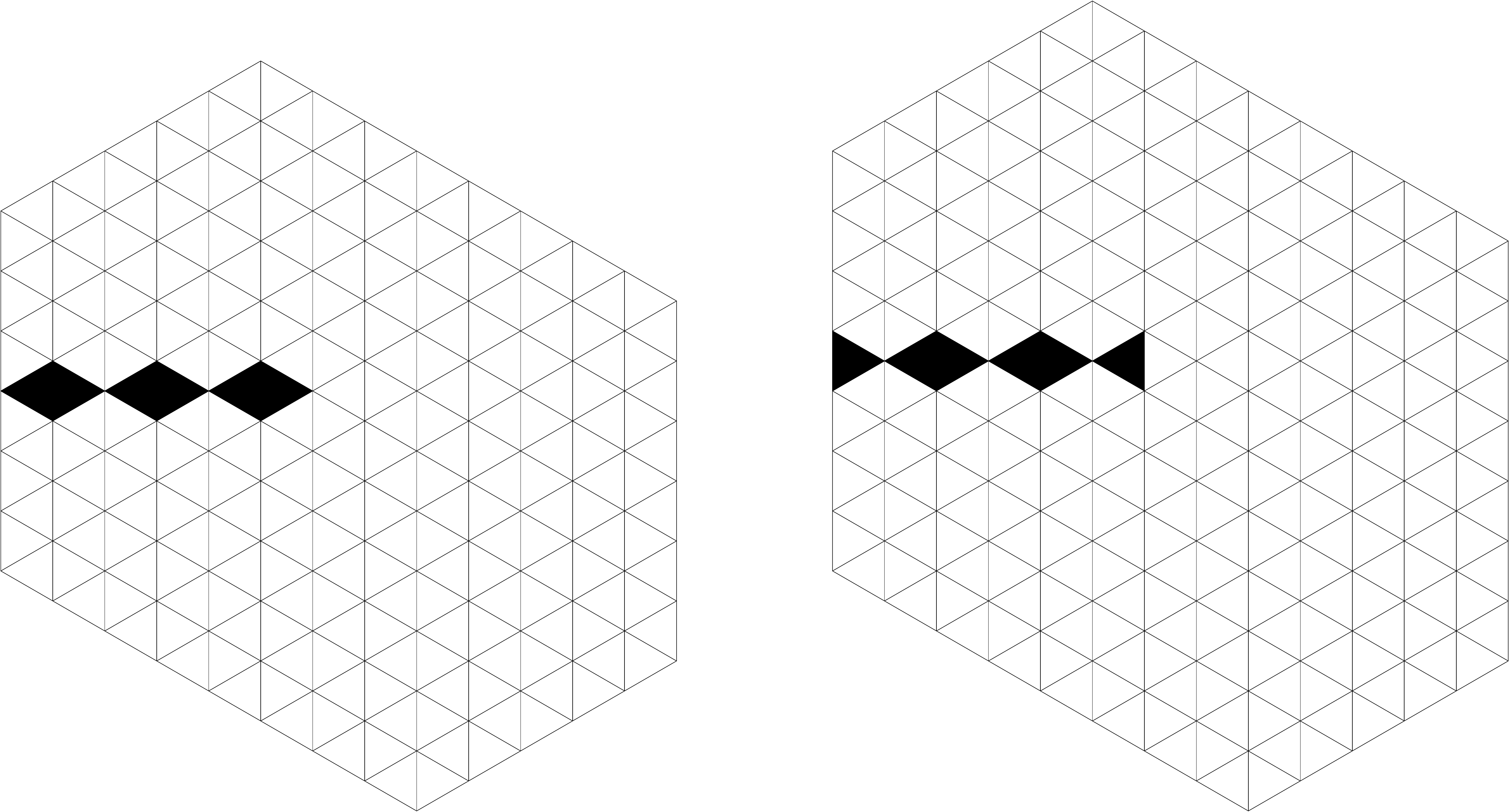}
    \caption{Two regions $H_{6,5,8;3}$ (left) and $H_{7,5,8;3}$ (right).}
\end{figure}

\begin{figure}
    \centering
    \includegraphics[width=16cm]{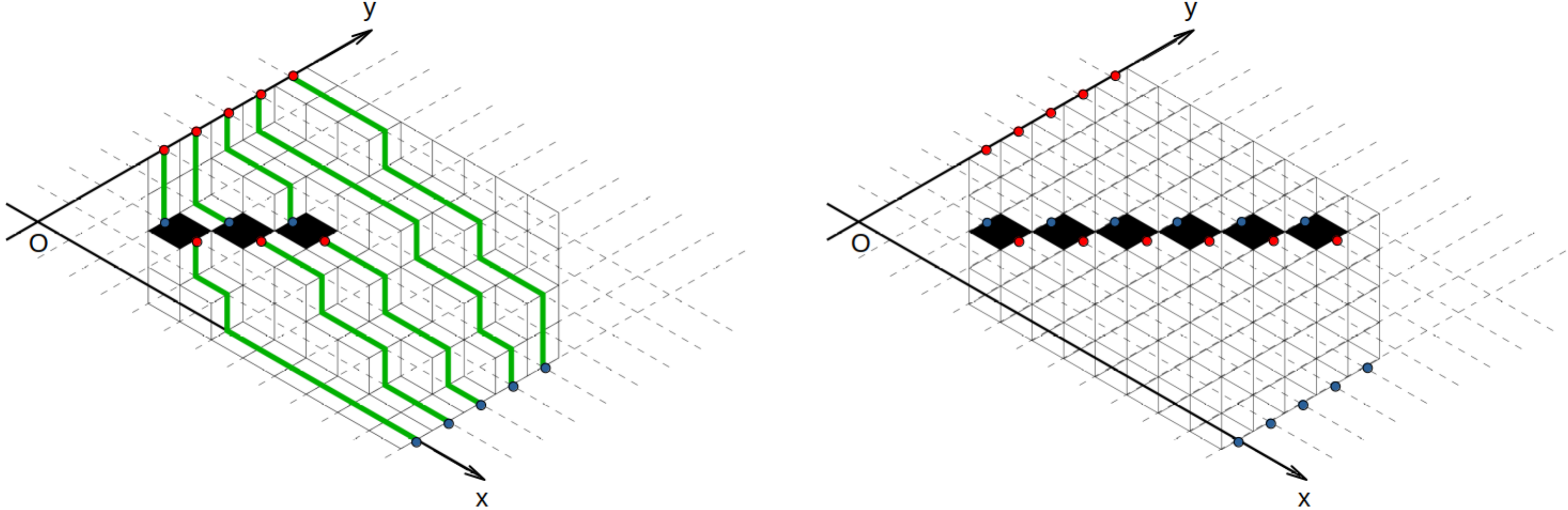}
    \caption{The left picture shows a lozenge tiling of the region $H_{4,5,8;3}$ and the corresponding family of nonintersecting lattice paths, where start and end points are marked by red and blue dots, respectively. On the right picture, one can easily see that the region $H_{4,5,8;6}$ does not allow such a family of nonintersecting paths that connect all start and end points and thus does not have any lozenge tilings.}
\end{figure}

To state the main theorem, we extend\footnote{This extension is based on the fact that the shifted factorial $(a)_n$ and the quotient of two gamma functions $\frac{\Gamma(a+n)}{\Gamma(a)}$ agree for a positive integer (or half-integer) $a$ and a nonnegative integer $n$. One can check that the definition of $(a)_{n}$ for a negative integer $n$ (=$\frac{1}{(a-1)(a-2)\cdots(a+n)}$) also coincides with the value of $\frac{\Gamma(a+n)}{\Gamma(a)}$.} the notion of shifted factorial mentioned earlier. For a positive integer (or half-integer) $a$ and an integer $n$ such that $a+n>0$, $(a)_n$ is defined as follows:
\begin{equation*}
    (a)_n:=
    \begin{cases}
        a(a+1)\cdots(a+n-1)                            & \text{if $n>0$,}\\
        1                            & \text{if $n=0$,}\\        
        \displaystyle\frac{1}{(a-1)(a-2)\cdots(a+n)}      & \text{if $n<0$ and $a+n>0$.}
    \end{cases}
\end{equation*}

Also, throughout this paper, an empty product (for example, $\displaystyle \prod_{k=0}^{n}[\cdots]$ with $n=-1$) is understood as $1$.

For any region $R$ on a triangular lattice, let $M(R)$ be the number of lozenge tilings of the region. If the region $R$ has empty interior, then we set $M(R)=1$ (in this case $R$ has one tiling -- the \textit{empty tiling}). The main theorem of this paper is as follows:

\begin{thm}
If $a$, $b$, $c$ and $d$ are nonnegative integers such that $d \leq b \leq c$ holds, then

\begin{equation}
    \frac{M(H_{2a,b,c;d})}{M(H_{2a,b,c})}=\prod_{k=0}^{d-1} \frac{(k+\frac{1}{2})_{b-2k}(a+k+1)_{b-2k-1}(b-k+\frac{1}{2})_{\lfloor\frac{c-b}{2}\rfloor}(c-k)_{-\lfloor\frac{c-b}{2}\rfloor}}{(a+k+\frac{1}{2})_{b-2k}(k+1)_{b-2k-1}(a+b-k+\frac{1}{2})_{\lfloor\frac{c-b}{2}\rfloor}(a+c-k)_{-\lfloor\frac{c-b}{2}\rfloor}}
\end{equation}
and
\begin{equation}
    \frac{M(H_{2a+1,b,c;d})}{M(H_{2a+1,b,c})}=\frac{1}{4^d}\prod_{k=0}^{d-1} \frac{(a+k+1)_{c-2k}(k+\frac{3}{2})_{c-2k-2}(b-k)_{\lfloor\frac{c-b}{2}\rfloor}(c-k-\frac{1}{2})_{-\lfloor\frac{c-b}{2}\rfloor}}{(k+1)_{c-2k-1}(a+k+\frac{3}{2})_{c-2k-1}(a+b-k+1)_{\lfloor\frac{c-b}{2}\rfloor}(a+c-k+\frac{1}{2})_{-\lfloor\frac{c-b}{2}\rfloor}}.
\end{equation}
\end{thm}

From MacMahon's theorem\footnote{It is $\displaystyle M(H_{a,b,c})=\prod_{i=1}^{a}\prod_{j=1}^{b}\prod_{k=1}^{c}\frac{i+j+k-1}{i+j+k-2}=\frac{H(a)H(b)H(c)H(a+b+c)}{H(a+b)H(b+c)H(c+a)}$,
where $H(n):=\displaystyle \prod_{i=0}^{n-1} i!$.}, we know the product formulas for $M(H_{2a,b,c})$ and $M(H_{2a+1,b,c})$. Hence, by multiplying $M(H_{2a,b,c})$ and $M(H_{2a+1,b,c})$ to (2.1) and (2.2), respectively, our theorem provides the product formulas for the numbers of lozenge tilings of the two regions $H_{2a,b,c;d}$ and $H_{2a+1,b,c;d}$. Thus, equations (2.1) and (2.2) provide a generalization of MacMahon's theorem since we recover MacMahon's theorem by replacing $d$ by $0$ in the formulas for $M(H_{2a,b,c;d})$ and $M(H_{2a+1,b,c;d})$ (see also [2], [3], [4], [7], [8], [19], [20], [21], [22], [23], [28] for similar generalizations in the literature).

\begin{figure}
    \centering
    \includegraphics[width=8cm]{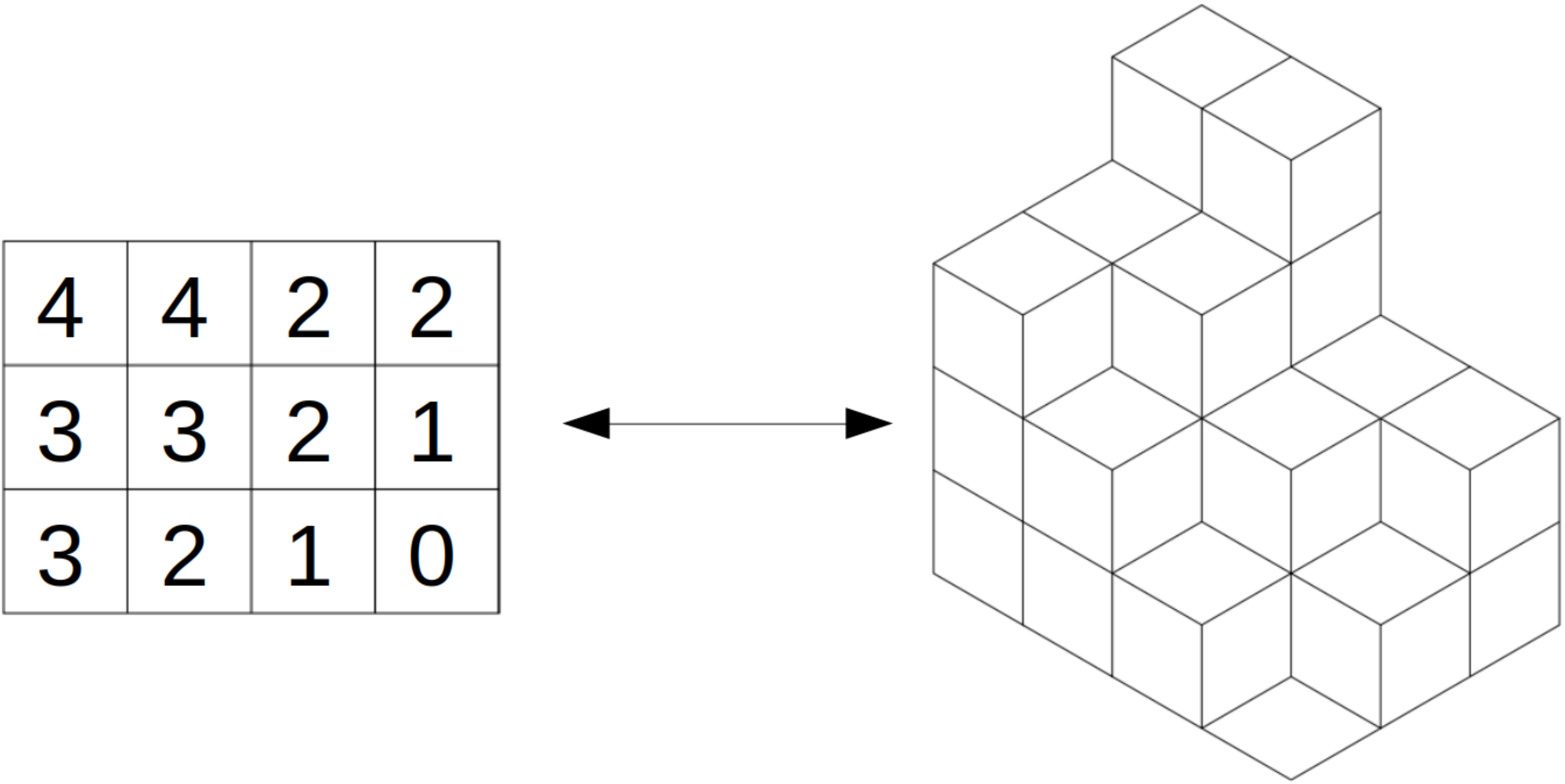}
    \caption{The left picture is an example of a plane partition, and the right picture is the 3-D interpretation of the plane partition as a stack of unit cubes.}
\end{figure}

\begin{figure}
    \centering
    \includegraphics[width=10cm]{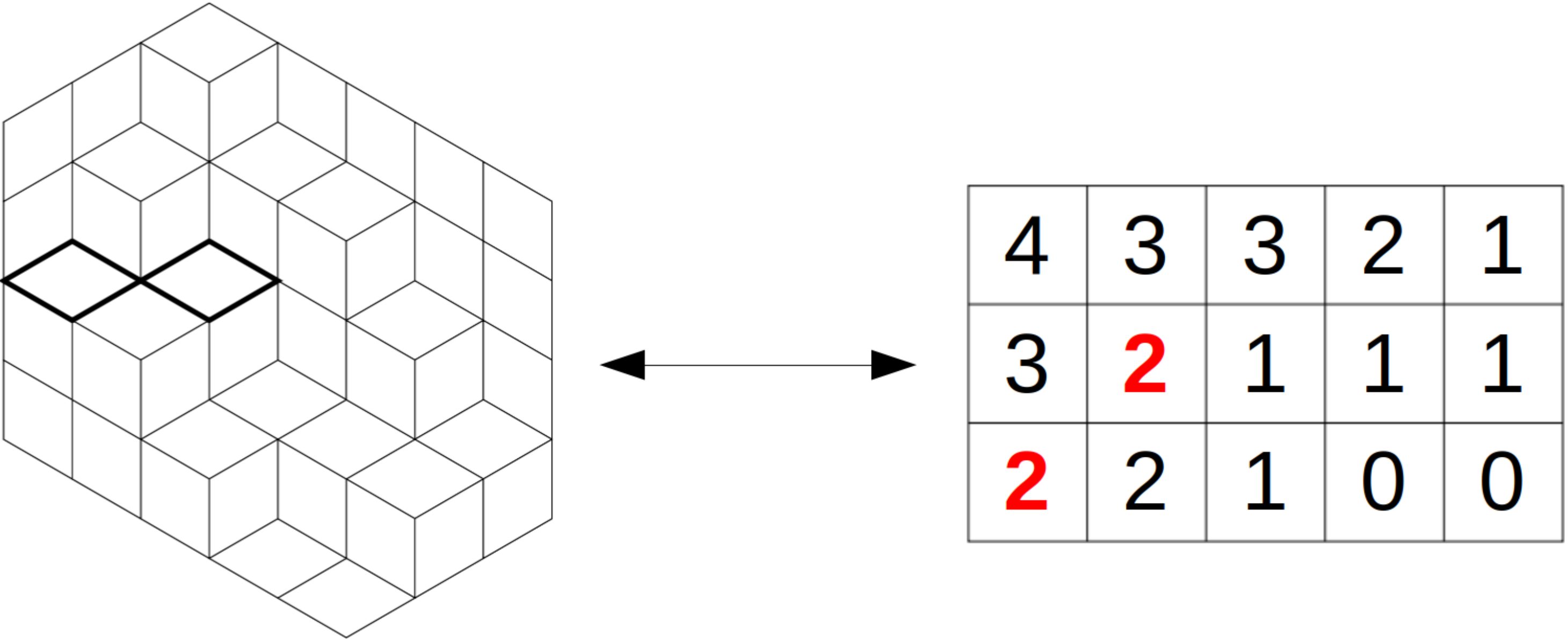}
    \caption{The figure shows the correspondence between lozenge tilings of $H_{4,3,5;2}$ (left) and plane partitions in $P(3,5,4;2)$ (right). We consider the lozenges with bolded edges either being removed or being fixed, as needed.}
\end{figure}

Equation (2.1) in the main theorem can also be interpreted in terms of plane partitions.
Recall that a \textit{plane partition} is an array of nonnegative integers $(\pi_{i,j})_{i,j\geq1}$ that is weakly decreasing along rows and down columns. More precisely, it satisfies (1) $\pi_{i,j}\geq \pi_{i+1,j}\,\, \forall i, j$, (2) $\pi_{i,j}\geq \pi_{i,j+1}\,\, \forall i, j$, and (3) $\pi_{i,j}=0$ for all but finitely many pairs $(i,j)$.
One typical interpretation of a plane partition is to view it as a stack of unit cubes. This can be obtained by stacking $n$ unit cubes on the position where $n$ is written in the plane partition (see Figure 2.3). From this visual interpretation, it is natural to consider a boxed plane partition. An \textit{$(a,b,c)$-boxed plane partition} is a plane partition with at most $a$ rows, at most $b$ columns, with entries that are at most $c$. Hence, the set of $(a,b,c)$-boxed plane partitions, denoted by $P(a,b,c)$, is $\{(\pi_{i,j})_{1\leq i\leq a, 1\leq j\leq b}|\pi_{i,j}\in\mathbb{Z}_{\geq0},\pi_{i,j}\geq \pi_{i+1,j}, \pi_{i,j}\geq \pi_{i,j+1}, c\geq\pi_{i,j} \,\, \forall i, j\}$.

Now, for nonnegative integers $a$, $b$, and $c$ such that $b \leq c$ holds, we consider $P(b,c,2a)$ and a certain restriction on it. For a nonnegative integer $d$ such that $d \leq b$ holds, let $P(b,c,2a;d)$ be the set of plane partitions in $P(b,c,2a)$ that satisfy the following additional conditions:
\begin{equation*}
    \pi_{b+1-i,i}=a \text{ for } i=1,2,\ldots,d.
\end{equation*}

One can easily see that there is a natural bijection between the set of lozenge tilings of the region $H_{2a,b,c;d}$ and the set of lozenge tilings of the region $H_{2a,b,c}$ with $d$ fixed lozenges on the perpendicular bisector of the left side, where positions of the fixed lozenges correspond to the positions of $d$ removed lozenges (or $2d$ triangles) on $H_{2a,b,c;d}$. Also, the bijection of David and Tomei gives a natural one-to-one correspondence between lozenge tilings of $H_{2a,b,c}$ with the $d$ fixed lozenges and plane partitions in $P(b,c,2a;d)$ (see Figure 2.4 that illustrates these bijections). Hence, the ratio $\displaystyle \frac{M(H_{2a,b,c;d})}{M(H_{2a,b,c})}$ is the same as $\displaystyle \frac{|P(b,c,2a;d)|}{|P(b,c,2a)|}$, and we obtain the following direct corollary of Theorem 2.1.

\begin{cor}
If $a$, $b$, $c$, and $d$ are nonnegative integers such that $d \leq b \leq c$ holds, then
\begin{equation}
    \frac{|P(b,c,2a;d)|}{|P(b,c,2a)|}=\prod_{k=0}^{d-1} \frac{(k+\frac{1}{2})_{b-2k}(a+k+1)_{b-2k-1}(b-k+\frac{1}{2})_{\lfloor\frac{c-b}{2}\rfloor}(c-k)_{-\lfloor\frac{c-b}{2}\rfloor}}{(a+k+\frac{1}{2})_{b-2k}(k+1)_{b-2k-1}(a+b-k+\frac{1}{2})_{\lfloor\frac{c-b}{2}\rfloor}(a+c-k)_{-\lfloor\frac{c-b}{2}\rfloor}}.
\end{equation}
\end{cor}

The explicit formulas in (2.1) and (2.2) motivate us to analyze the limit or asymptotic behavior of them (the same kind of analysis was already carried out in the literature. For example, see [5], [10], [12], and [13]). To consider the scaling limit, we scale the parameters $a, b, c$, and $d$ by a factor $N$ and observe the limit. Thus, for a positive integer $N$, we consider the ratios $\displaystyle \frac{M(H_{2aN,bN,cN;dN})}{M(H_{2aN,bN,cN})}$ and $\displaystyle \frac{M(H_{2aN+1,bN,cN;dN})}{M(H_{2aN+1,bN,cN})}$.

Numerical data suggests that these ratios \textit{usually} (but not always!) tend to $0$ as $N$ approaches infinity. Although they do not have interesting limits, one can then ask about their asymptotic behavior. Using Stirling's formula\footnote{$n! \sim \sqrt{2\pi n}(\frac{n}{e})^n$ as $n\to\infty$.} and the definition of the Glaisher-Kinkelin constant (see [15]), we can show the following results on the behavior of these two ratios for large $N$. To state the theorem, we first recall the definitions of the Barnes $G$-function and the Glaisher-Kinkelin constant $A$.
For an integer $n$, the value of the Barnes $G$-function is given as follows:
\begin{equation}
    G(n):=
\begin{cases}
    0                                                    &  \text{if $n$ is negative or $0$,}\\
    \displaystyle \prod_{i=0}^{n-2}i!                    &  \text{if $n$ is positive.}
\end{cases}
\end{equation}
The Glaisher-Kinkelin constant $A$ is defined using the Barnes $G$-function as follows:
\begin{equation}
    A:=\lim_{n\to\infty}\frac{(2\pi)^{\frac{n}{2}}n^{\frac{1}{2}n^2-\frac{1}{12}}e^{-\frac{3}{4}n^2+\frac{1}{12}}}{G(n+1)}=1.282427129\cdots.
\end{equation}
We can now state the following theorem on the asymptotic behavior of the two ratios above for large $N$. To simplify the analysis, we additionally assume that $b$ and $c$ have the same parity.

\begin{thm}
Let $a$, $b$, $c$, and $d$ be positive integers such that 1) $d < b \leq c$ holds and 2) $b$ and $c$ have the same parity. As $N$ grows large,\footnote{Throughout this paper, $a(N)\sim b(N)$ (as $N\to\infty$) means $\displaystyle \lim_{N\to\infty}\frac{a(N)}{b(N)}=1$.}
\begin{equation}
    \frac{M(H_{2aN,bN,cN;dN})}{M(H_{2aN,bN,cN})} \sim \frac{2^{\frac{7}{24}}e^{\frac{1}{8}}}{A^{\frac{3}{2}}N^{\frac{1}{8}}} K_{1}(a,b,c,d) K_{2}^{\frac{N}{2}}(a,b,c,d) K_{3}^{\frac{N^{2}}{2}}(a,b,c,d)\\
\end{equation}
and
\begin{equation}
    \frac{M(H_{2aN+1,bN,cN;dN})}{M(H_{2aN+1,bN,cN})} \sim \frac{e^{\frac{1}{8}}}{2^{\frac{5}{24}}A^{\frac{3}{2}}N^{\frac{1}{8}}} L_{1}(a,b,c,d) L_{2}^{\frac{N}{2}}(a,b,c,d) K_{3}^{\frac{N^{2}}{2}}(a,b,c,d),\\
\end{equation}
where
\begin{equation}
    K_{1}(a,b,c,d)=\Bigg[\frac{a+d}{ad}\Bigg]^{\frac{1}{8}}\Bigg[\frac{bc(a+b-d)(a+c-d)}{(a+b)(a+c)(b-d)(c-d)}\Bigg]^{\frac{1}{12}}\Bigg[\frac{(a+\frac{1}{2}b+\frac{1}{2}c)(\frac{1}{2}b+\frac{1}{2}c-d)}{(\frac{1}{2}b+\frac{1}{2}c)(a+\frac{1}{2}b+\frac{1}{2}c-d)}\Bigg]^{\frac{1}{24}},
\end{equation}
\begin{equation}
    K_{2}(a,b,c,d)=\frac{a^{a}d^{d}(\frac{1}{2}b+\frac{1}{2}c)^{\frac{1}{2}b+\frac{1}{2}c}(a+\frac{1}{2}b+\frac{1}{2}c-d)^{a+\frac{1}{2}b+\frac{1}{2}c-d}}{(a+d)^{a+d}(a+\frac{1}{2}b+\frac{1}{2}c)^{a+\frac{1}{2}b+\frac{1}{2}c}(\frac{1}{2}b+\frac{1}{2}c-d)^{\frac{1}{2}b+\frac{1}{2}c-d}},
\end{equation}
\begin{equation}
\begin{aligned}
    K_{3}(a,b,c,d)=&\frac{(\frac{1}{2}b+\frac{1}{2}c)^{2(\frac{1}{2}b+\frac{1}{2}c)^2}}{b^{b^2}c^{c^2}} \frac{(a+\frac{1}{2}b+\frac{1}{2}c-d)^{2(a+\frac{1}{2}b+\frac{1}{2}c-d)^2}}{(a+b-d)^{(a+b-d)^2}(a+c-d)^{(a+c-d)^2}}\\ &\times\frac{(b-d)^{(b-d)^2}(c-d)^{(c-d)^2}}{(\frac{1}{2}b+\frac{1}{2}c-d)^{2(\frac{1}{2}b+\frac{1}{2}c-d)^2}} \frac{(a+b)^{(a+b)^2}(a+c)^{(a+c)^2}}{(a+\frac{1}{2}b+\frac{1}{2}c)^{2(a+\frac{1}{2}b+\frac{1}{2}c)^2}},    
\end{aligned}
\end{equation}
\begin{equation}
\begin{aligned}
    &L_{1}(a,b,c,d)\\
    &=\Bigg[\frac{a+\frac{1}{2}b+\frac{1}{2}c-d}{a+\frac{1}{2}b+\frac{1}{2}c}\Bigg]^{\frac{11}{24}}\Bigg[\frac{(a+b)(a+c)}{(a+b-d)(a+c-d)}\Bigg]^{\frac{5}{12}}\Bigg[\frac{a+d}{ad}\Bigg]^{\frac{1}{8}}\Bigg[\frac{bc}{(b-d)(c-d)}\Bigg]^{\frac{1}{12}}\Bigg[\frac{\frac{1}{2}b+\frac{1}{2}c-d}{\frac{1}{2}b+\frac{1}{2}c}\Bigg]^{\frac{1}{24}},
\end{aligned}
\end{equation}
and
\begin{equation}
\begin{aligned}
    L_{2}(a,b,c,d)=&
    \Bigg[\frac{(a+d)^{a+d}}{a^{a}d^{d}}\frac{(a+\frac{1}{2}b+\frac{1}{2}c)^{a+\frac{1}{2}b+\frac{1}{2}c}(\frac{1}{2}b+\frac{1}{2}c-d)^{\frac{1}{2}b+\frac{1}{2}c-d}}{(a+\frac{1}{2}b+\frac{1}{2}c-d)^{a+\frac{1}{2}b+\frac{1}{2}c-d}(\frac{1}{2}b+\frac{1}{2}c)^{\frac{1}{2}b+\frac{1}{2}c}}\Bigg]\\
    &\times\Bigg[\frac{1}{4^{d}}\frac{(a+b)^{a+b}(a+c)^{a+c}}{(a+\frac{1}{2}b+\frac{1}{2}c)^{2(a+\frac{1}{2}b+\frac{1}{2}c)}}\frac{(a+\frac{1}{2}b+\frac{1}{2}c-d)^{2(a+\frac{1}{2}b+\frac{1}{2}c-d)}}{(a+b-d)^{a+b-d}(a+c-d)^{a+c-d}}\Bigg]^{2}.
\end{aligned}
\end{equation}
\end{thm}

We postpone the proof of Theorem 2.3 to the last section. There are two simple facts about these formulas. One is that the five constants $K_{1}(a,b,c,d), K_{2}(a,b,c,d), K_{3}(a,b,c,d), L_{1}(a,b,c,d)$, and $L_{2}(a,b,c,d)$ are symmetric in the variables $b$ and $c$. Another is that $K_{3}(a,b,c,d)$ becomes $1$ when $b$ and $c$ are the same. Thus, the two ratios are simplified as follows when $b=c$:
\begin{equation}
    \frac{M(H_{2aN,cN,cN;dN})}{M(H_{2aN,cN,cN})} \sim \frac{2^{\frac{7}{24}}e^{\frac{1}{8}}}{A^{\frac{3}{2}}N^{\frac{1}{8}}} \Bigg[\frac{(a+d)c(a+c-d)}{ad(a+c)(c-d)}\Bigg]^{\frac{1}{8}} \Bigg[\frac{a^{a}d^{d}c^{c}(a+c-d)^{a+c-d}}{(a+d)^{a+d}(a+c)^{a+c}(c-d)^{c-d}}\Bigg]^{\frac{N}{2}}\\
\end{equation}
and
\begin{equation}
\begin{aligned}
    \frac{M(H_{2aN+1,cN,cN;dN})}{M(H_{2aN+1,cN,cN})} \sim \frac{e^{\frac{1}{8}}}{2^{\frac{5}{24}}A^{\frac{3}{2}}N^{\frac{1}{8}}} \Bigg[\frac{(a+d)c}{ad(c-d)}\Bigg]^{\frac{1}{8}} \Bigg[\frac{a+c}{a+c-d}\Bigg]^{\frac{3}{8}}\Bigg[\frac{(a+d)^{a+d}(a+c)^{a+c}(c-d)^{c-d}}{4^{2d}a^{a}d^{d}(a+c-d)^{a+c-d}c^{c}}\Bigg]^{\frac{N}{2}}.
\end{aligned}
\end{equation}

As a direct corollary of Theorem 2.3, we can find the limit of the ratios as $N$ approaches infinity. The proof easily follows from the two inequalities $K_{3}(a,b,c,d)\leq 1$ (the inequality is strict unless $b=c$) and $K_{2}(a,c,c,d)< 1$, which can be proved using convexity of the functions $\displaystyle f(t)=ln\Bigg[\frac{(\frac{1}{2}b+\frac{1}{2}c+t)^{2(\frac{1}{2}b+\frac{1}{2}c+t)^2}}{(b+t)^{(b+t)^2}(c+t)^{(c+t)^2}}$\Bigg] and $g(t)=ln(t^t)=t\,ln\,t$, respectively\footnote{For example, if one take natural logarithm to the latter inequality $K_{2}(a,c,c,d)< 1$, then one can see that it is enough to show that $g(a+d)>g(a)+g(d)$ and $g(a+c)+g(c-d)>g(a+c-d)+g(c)$ hold. These two inequalities can be verified using Jensen's inequality and the fact that $\displaystyle\lim_{x\rightarrow 0^{+}}g(x)=0$.}.

\begin{cor}
If $a$, $b$, $c$, and $d$ are positive integers such that 1) $d < b \leq c$ holds and 2) $b$ and $c$ have the same parity, then
\begin{equation}
    \lim_{N\rightarrow\infty}\frac{M(H_{2aN,bN,cN;dN})}{M(H_{2aN,bN,cN})} = 0\\
\end{equation}
and
\begin{equation}
    \lim_{N\rightarrow\infty}\frac{M(H_{2aN+1,bN,cN;dN})}{M(H_{2aN+1,bN,cN})}=
\begin{cases}
    \infty       &  \text{if $b=c$ and $\displaystyle \frac{(a+d)^{a+d}(a+c)^{a+c}(c-d)^{c-d}}{4^{2d}a^{a}d^{d}(a+c-d)^{a+c-d}c^{c}} > 1$}, \\
    0            &  \text{otherwise.}
\end{cases}
\end{equation}
\end{cor}
Hence, the behaviors of the two ratios are completely different when $a, b, c,$ and $d$ satisfy the conditions described on the right side of (2.16). Note that there exist 4-tuples $(a,b,c,d)$ that satisfy the conditions. For an example, $(a,b,c,d)=(3,3,3,1)$ satisfies them.

\section{Preparation for the Proof}

In this section, we list the theorems that we use to prove the main theorem, Theorem 2.1. 

To state the first theorem, we define the two regions $R_{m,n,x}$ and $\overline{R}_{m,n,x}$ on a triangular lattice (see Figure 3.1 and Figure 3.2).

We first define the region $R_{m,n,x}$. Let \emph{O} be any lattice point on the triangular lattice and \emph{$\Bar{O}$} be the lattice point that is one unit northwest from \emph{O}. Consider the horizontal line \emph{l} through \emph{O} (which is not a lattice line) and let \emph{A} be the $n$th lattice point to the right of \emph{O} which lies on \emph{l} (when $n=0$, \emph{A}=\emph{O}). Similarly, consider the horizontal line $\Bar{\emph{l}}$ through \emph{$\Bar{O}$} and let \emph{B} be the $m$th lattice point left of \emph{$\Bar{O}$} (when $m=0$, \emph{B}=\emph{$\Bar{O}$}). The region $R_{m,n,x}$ is defined differently, depending on whether $m=0$ or $m\geq1$.

\begin{figure}
    \centering
    \includegraphics[width=15cm]{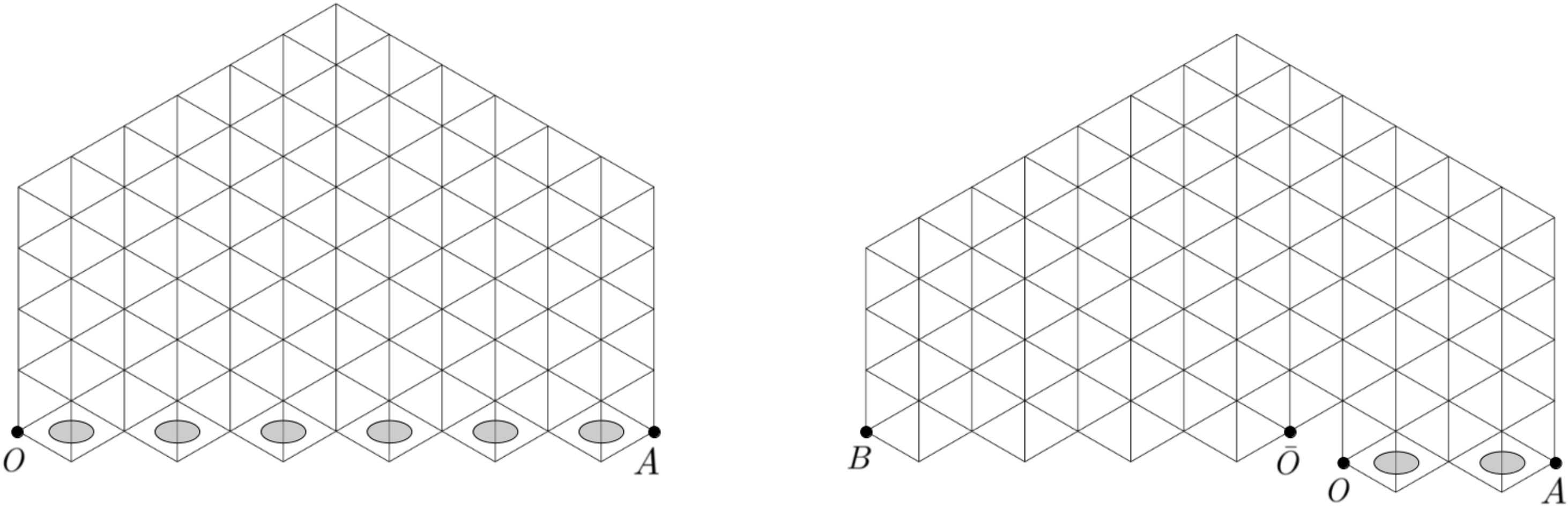}
    \caption{The regions $R_{0,6,3}$ (left) and $R_{4,2,3}$ (right). The lozenges that lie in the position indicated by the shaded ellipses are given a weight of $\frac{1}{2}$, and all other lozenges have weight 1.}
\end{figure}

When $m=0$, the region $R_{0,n,x}$ is defined as follows without using \emph{$\Bar{O}$}. From \emph{A} to \emph{O}, we follow the zigzag line along \emph{l}, by alternating moving one unit to the southwest and one unit to the northwest. Next, we move $x+1$ units to the north, $n$ units to the northeast, $n$ units to the southeast, and then $x+1$ units to the south until we reach \emph{A}. For any integer $x$ such that $x\geq-1$, $R_{0,n,x}$ is defined to be the bounded region enclosed by the path described above. In any tiling of this region, lozenges that lie along the line \emph{l} have weight $\frac{1}{2}$, while all other lozenges have weight 1 (see the left picture in Figure 3.1).

When $m\geq1$, from \emph{A} to \emph{O}, we follow the zigzag line along \emph{l}, as we did in the case when $m=0$. Next, from \emph{O} to \emph{$\Bar{O}$}, we move one unit to the north and then one unit to the southwest. Now, we connect \emph{$\Bar{O}$} and \emph{B} using the same type of zigzag line along $\Bar{\emph{l}}$. From \emph{B}, we move $x$ units to the north, $m+n+1$ units to the northeast, $m+n$ units to the southeast, and then $x+1$ units to the south until we reach \emph{A}. In this case, for any nonnegative integer $x$, $R_{m,n,x}$ is defined to be the bounded region enclosed by the path we just described. Like the previous case, in any tilings of this region, lozenges that lie along the line \emph{l} are given a weight $\frac{1}{2}$ (all other lozenges have weight 1. See the right picture in Figure 3.1).

\begin{figure}
    \centering
    \includegraphics[width=15cm]{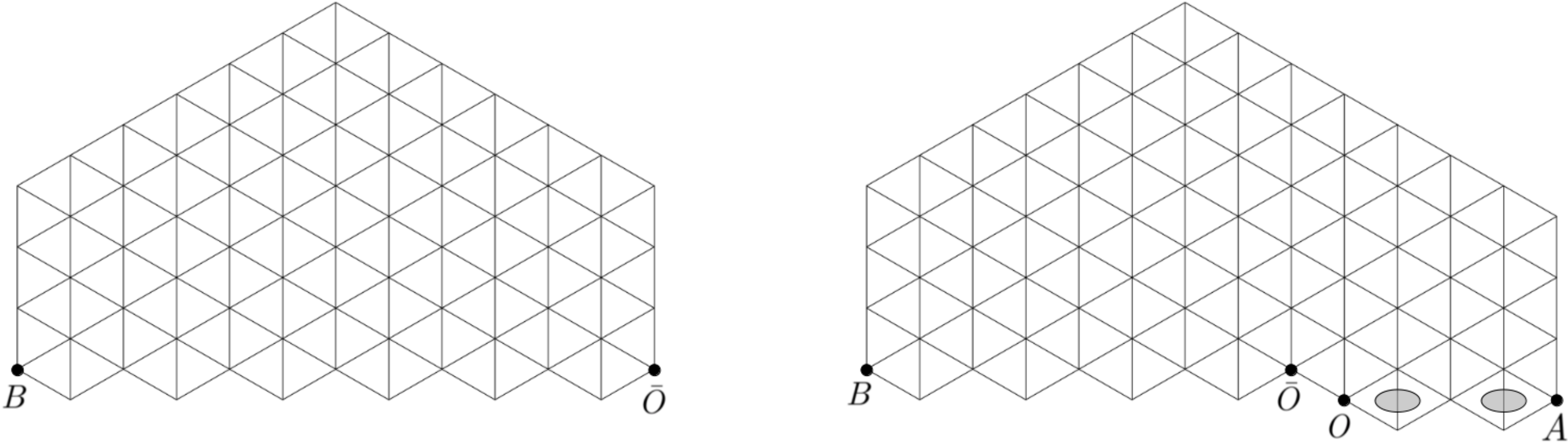}
    \caption{The regions $\overline{R}_{6,0,3}$ (left) and $\overline{R}_{4,2,3}$ (right). The lozenges that lie in the position indicated by the shaded ellipses are given a weight of $\frac{1}{2}$, and all other lozenges have weight 1.}
\end{figure}

The region $\overline{R}_{m,n,x}$ is defined similarly (see Figure 3.2). The differences are
\begin{enumerate}
    \item {} we connect \emph{O} and \emph{$\Bar{O}$} using the unit segment connecting them.
    \item {} when $n=0$, from \emph{B}, we move $x$ units to the north, $m$ units to the northeast, $m$ units to the southeast, and then $x$ units to the south until we reach \emph{$\Bar{O}$} (so \emph{O} is not used when we define the region $\overline{R}_{m,0,x}$). In this case, no lozenges lie along the line \emph{l}, so every lozenge have weight 1 (see the left picture in Figure 3.2).
    \item {} when $n\geq1$, from \emph{B}, we move $x$ units to the north, $m+n$ units to the northeast, $m+n+1$ units to the southeast, and then $x$ units to the south until we reach \emph{A}. As before, lozenges that lie along the line \emph{l} have weight $\frac{1}{2}$ and all other lozenges have weight 1 (see the right picture in Figure 3.2).
\end{enumerate}

Tilings of two regions $R_{m,n,x}$ and $\overline{R}_{m,n,x}$ may well contain lozenges that have weight $\frac{1}{2}$. In general, when we assign weights to lozenges in a tiling of a region R, we use $M(R)$ to denote the sum over weighted tilings of $R$, where the weight of a tiling is the product of the weights of the tiles that comprise it. Note that if all lozenges in the region have weight $1$, then every tiling is weighted by $1$ and $M(R)$ simply enumerates tilings of $R$.

The following is a special case of Proposition 2.1 in [2]. In [2], Ciucu provided product formulas for more general regions (see also [24], where Lai and Rohatgi provided elegant weighted generalizations of the result of Ciucu [2]). To state the special case of the formulas of Ciucu, we need to define a notation $<a, a+n>$ for a positive integer (or half-integer) $a$ and an integer $n$.
For a positive integer (or half-integer) $a$ and an integer $n$, let

\begin{equation}
\begin{aligned}
    <a, a+n>&:=
\begin{cases}
    \displaystyle \prod_{i=0}^{n}(a+i)^{\min(i+1, n+1-i)} & \text{if $n$ is positive or $0$,}\\
    1                                                    &  \text{if $n$ is negative.}
\end{cases}
\\
    &=a(a+1)^2(a+2)^3\cdots(a+n-2)^3(a+n-1)^2(a+n).
\end{aligned}
\end{equation}

\begin{thm}[{[2]}, Proposition 2.1] For nonnegative integers $m, n,$ and $x$,

\begin{equation}
\begin{aligned}
    &M(R_{m,n,x})\\
    &=2^{n(n-1)/2-2mn}\prod_{i=1}^{m}\frac{1}{(2i)!}\prod_{i=1}^{n}\frac{1}{(2i-1)!}\frac{\displaystyle \prod_{1\leq i< j\leq m}(j-i)\prod_{1\leq i< j\leq n}(j-i)}{\displaystyle \prod_{i=1}^{m}\prod_{j=1}^{n}(i+j)}\\
    &\times(x+n+1)_{m}(x+n+2)_{m}<x+2,x+n><x+\frac{3}{2},x+\frac{2n+1}{2}>\\
    &\times\prod_{i=1}^{n}\frac{(x+i)_{m}}{(x+i+\frac{1}{2})_{m}}\prod_{i=1}^{m}(2x+n+i+2)_{n+i-1}
\end{aligned}
\end{equation}
and
\begin{equation}
\begin{aligned}
    &M(\overline{R}_{m,n,x})\\
    &=2^{m(m-1)/2-2mn-n}\prod_{i=1}^{m}\frac{1}{(2i-1)!}\prod_{i=1}^{n}\frac{1}{(2i)!}\frac{\displaystyle \prod_{1\leq i< j\leq m}(j-i)\prod_{1\leq i< j\leq n}(j-i)}{\displaystyle \prod_{i=1}^{m}\prod_{j=1}^{n}(i+j)}\\
    &\times(x+m+1)_{n}<x+1,x+m><x+\frac{3}{2},x+\frac{2m-1}{2}>\\
    &\times\prod_{i=1}^{m}\frac{(x+i)_{n}}{(x+i+\frac{1}{2})_{n}}\prod_{i=1}^{n}(2x+m+i+1)_{m+i}.
\end{aligned}
\end{equation}
The equation $(3.2)$ is still valid when $m=0$ and $x=-1$.
\end{thm}

\textbf{Remark.} As we mentioned earlier, Theorem 3.1 is a special case of the more general result of Ciucu. The products over partitions $\lambda$ and $\mu$ in the paper of Ciucu [2] do not appear in our formulas because both $\lambda$ and $\mu$ are empty partitions in our case.\\

For our proof of Theorem 2.1, we will also require the Matching Factorization Theorem of Ciucu [1]. Instead of giving the precise statement of the theorem, we present how the theorem is applied to the regions $H_{2a,c,c;d}$ and $H_{2a+1,c,c;d}$ (see [1] for the precise statement of the Factorization Theorem).
We first split the region $H_{2a,c,c;d}$ into two subregions as follows: when $d=c$, the region consists of two components (this is because removed triangles connect the middle points of the left and the right sides). Let $H_{2a,c,c;c}^-$ be the region on the top and $H_{2a,c,c;d}^+$ be the region on the bottom. When $0\leq d<c$, from the midpoint of the right side of the region, we consider a zigzag lattice path of length $2c-2d$ whose direction alternates between the northwest and the southwest (see the pictures on the left in Figure 3.3). Now, we cut the region $H_{2a,c,c;d}$ along the zigzag lattice line. The region is now divided into two subregions. We denote the subregion on the top by $H_{2a,c,c;d}^-$. On the bottom subregion, we give weight $\frac{1}{2}$ to lozenges on the horizontal symmetry axis of $H_{2a,c,c;d}$ and denote it by $H_{2a,c,c;d}^+$.
Then, the following equation is true for any $d\in \{0, 1, \ldots, c\}$:
\begin{equation}
    M(H_{2a,c,c;d})=2^{c-d}M(H_{2a,c,c;d}^+)M(H_{2a,c,c;d}^-).
\end{equation}
For $d\in \{0, 1, \ldots, c-1\}$, it is due to the Factorization Theorem. For $d=c$, it comes from a simple fact that the number of lozenge tilings of a region with two components is the same as the product of the numbers of lozenges tilings of each component.

\begin{figure}
    \centering
    \includegraphics[width=10cm]{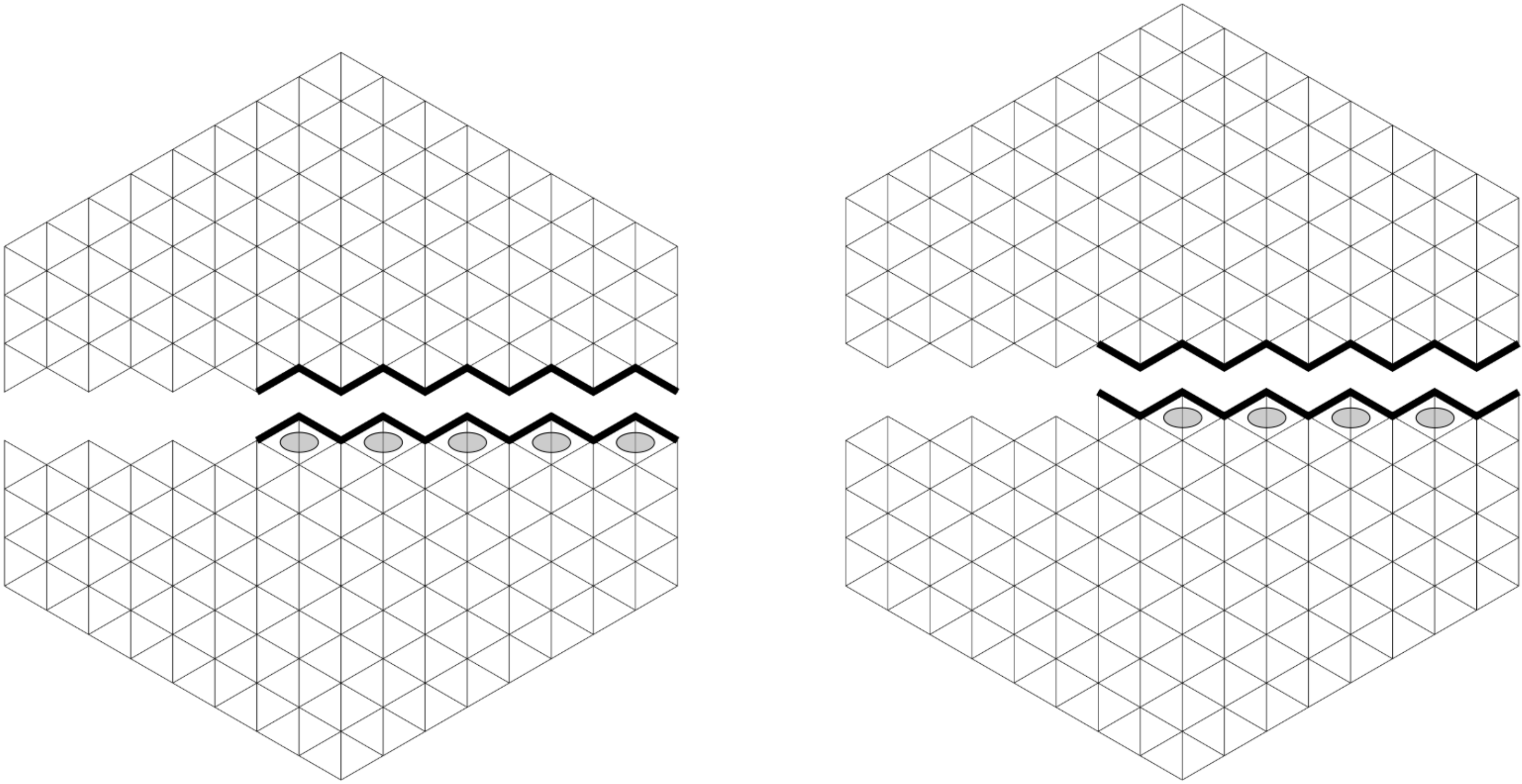}
    \caption{$H_{6,8,8;3}^+$ (bottom left), $H_{6,8,8;3}^-$ (top left), $H_{7,8,8;3}^+$ (bottom right), and $H_{7,8,8;3}^-$ (top right). The bold lines show how we cut the regions $H_{6,8,8;3}$ and $H_{7,8,8;3}$.}
\end{figure}

Similarly, for any $d\in \{0, 1, \ldots, c\}$, we split $H_{2a+1,c,c;d}$ into the two subregions $H_{2a+1,c,c;d}^+$ and $H_{2a+1,c,c;d}^-$ (see the pictures on the right in Figure 3.3). According to the Matching Factorization Theorem:

\begin{equation}
    M(H_{2a+1,c,c;d})=2^{c-d}M(H_{2a+1,c,c;d}^+)M(H_{2a+1,c,c;d}^-).
\end{equation}

The last theorems we present in this section are two versions of Kuo's graphical condensation [17]. As explained by Fulmek [11], this method is indeed an alternative combinatorial interpretation of the Kasteleyn-Percus method [16, 26]. We denote a bipartite graph by $G=(V_1, V_2, E)$, where $E$ is the set of edges of the graph and $(V_1, V_2)$ is the partition of the vertex set of the graph $G$ such that every edge in $E$ connects a vertex in $V_1$ and a vertex in $V_2$. Also, for any set of vertices $\{x_1, x_2,...,x_n\}$, let $G-\{x_1, x_2,...,x_n\}$ be the subgraph obtained by deleting $x_1, x_2,...,x_n$, together with all edges adjacent to those vertices. For any weighted graph $G$, let $M(G)$ denote the sum over weighted perfect matchings\footnote{\textit{A perfect matching} is a subset of edges of a bipartite graph such that every vertex is incident to precisely one edge. See the right picture in Figure 3.4 for an example.} of $G$, where the weight of a perfect matching is the product of the weights of the edges that comprise it. Observe that when all edges have weight $1$, $M(G)$ is simply the number of perfect matchings of $G$.

\begin{thm}[{[17]}, Theorem 2.1]
Let $G=(V_1, V_2, E)$ be a plane bipartite graph in which $|V_1|=|V_2|$. Let vertices $x, y, z,$ and $w$ appear in a cyclic order on a face of $G$. If $x,z \in V_1$ and $y,w \in V_2$, then
\begin{equation}
    M(G)M(G-\{x,y,z,w\})=M(G-\{x,y\})M(G-\{z,w\})+M(G-\{x,w\})M(G-\{y,z\}).
\end{equation}
\end{thm}

\begin{thm}[{[17]}, Theorem 2.4]
Let $G=(V_1, V_2, E)$ be a plane bipartite graph in which $|V_1|=|V_2|+1$. Let vertices $x, y, z,$ and $w$ appear in a cyclic order on a face of $G$. If $x,y,z \in V_1$ and $w \in V_2$, then
\begin{equation}
    M(G-\{y\})M(G-\{x,z,w\})=M(G-\{x\})M(G-\{y,z,w\})+M(G-\{z\})M(G-\{x,y,w\}).
\end{equation}
\end{thm}

\begin{figure}
    \centering
    \includegraphics[width=10cm]{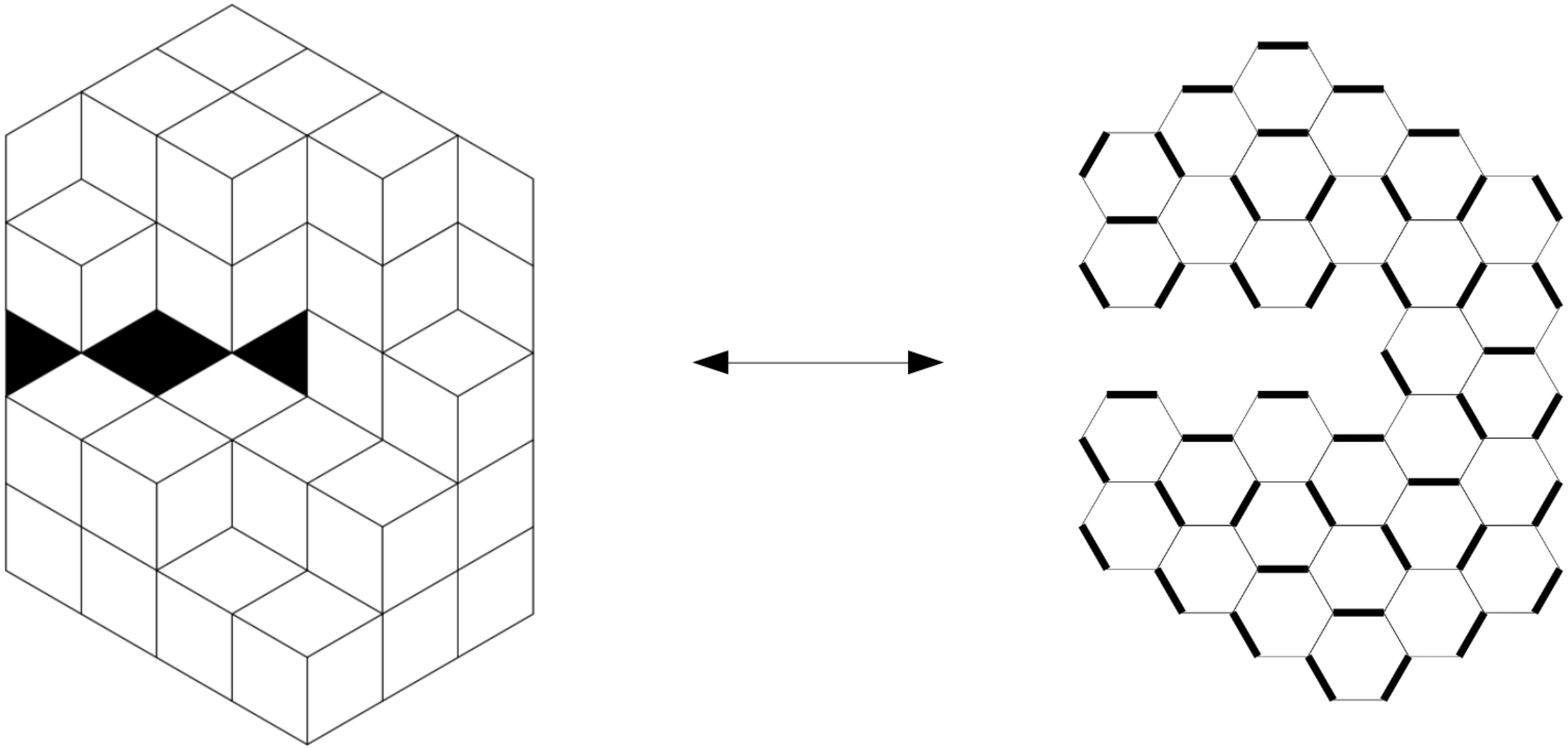}
    \caption{A lozenge tilling of the region $H_{5,3,4;2}$ (left) and the corresponding perfect matching of its dual graph (right).}
\end{figure}

A lozenge tiling of a region on a triangular lattice can be identified with a perfect matching of the dual graph of the region (on the hexagonal lattice). Given a lozenge tiling of the region, we choose every edge on the dual graph whose two vertices correspond to the adjacent unit triangles covered by a lozenge in the given tiling of the region (see Figure 3.4 that illustrates the correspondence. See also [18] for more details about it). This correspondence is bijective. Using the two theorems above and the bijection between perfect matchings and tilings, we construct two recurrence relations involving the number of lozenge tilings of the two regions $H_{2a,b,c;d}$ and $H_{2a+1,b,c;d}$. Next, using induction, we show that the numbers of lozenge tilings of these regions are given by (2.1) and (2.2). 

\begin{figure}
    \centering
    \includegraphics[width=10cm]{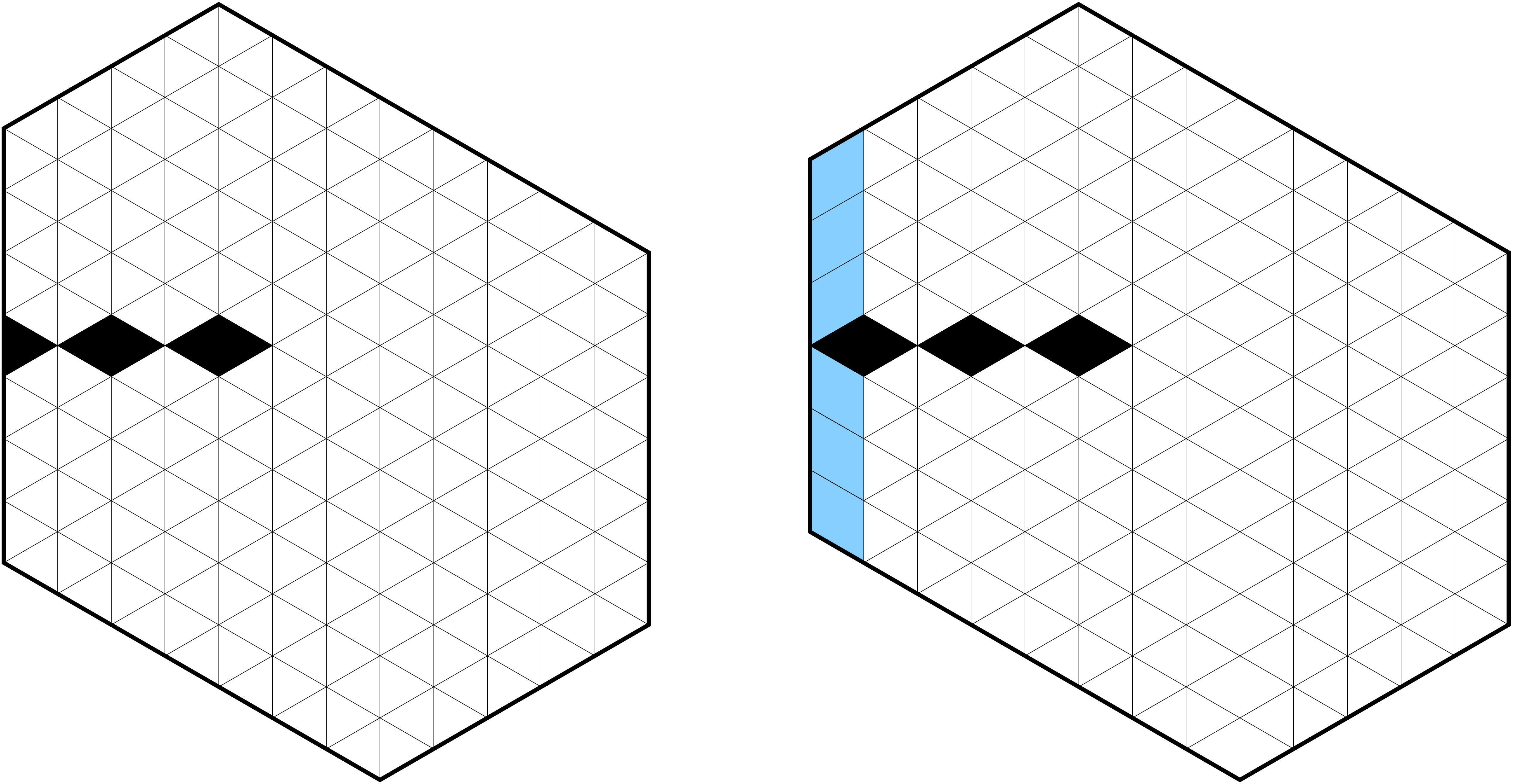}
    \caption{The regions with the same number of lozenge tilings. Forced lozenges are indicated by shading.}
\end{figure}

A lozenge-shaped tile on a region $R$ is a \textit{forced lozenge} if it is contained in every tiling of the region. In this paper, forced lozenges will be indicated by shading. If we denote the forced lozenge and its weight by $I$ and $\gamma$, respectively, since it is always part of the tilings, then one can easily see that $M(R-I)=M(R)/\gamma$ holds. In particular, if the forced lozenge has weight $\gamma=1$, then $M(R-I)=M(R)$. Figure 3.5 shows two regions that have the same number of tilings. They have the same number of tilings because the left region is obtained from the right one by deleting the forced lozenges. Since the right region is $H_{6,5,8;3}$, the number of lozenge tilings of the left region is $M(H_{6,5,8;3})$. In the proof of Theorem 2.1, we use this idea several times.

\section{A proof of Theorem 2.1}

The proof of Theorem 2.1 is organized as follows: first, using the idea provided by Ciucu and Krattenthaler in [6], we show the special case of Theorem 2.1 when $c-b=0$. Next, we construct two recurrence relations using Kuo's graphical condensation. By specializing in the recurrences, we give straightforward proof of the case when $c-b=1$. Lastly, we show the general case (when $c-b\geq 2$) using a double induction argument. In the proof, we assume $b>0$ since $b=0$ implies $d=0$, and we already know that equations (2.1) and (2.2) hold when $d=0$.\\

\textbf{Step 1: the case } $\mathbf{c-b=0}$

When $c-b=0$, (2.1) and (2.2) become the following identities:

\begin{equation}
    \frac{M(H_{2a,c,c;d})}{M(H_{2a,c,c})}= \prod_{k=0}^{d-1}  \frac{(k+\frac{1}{2})_{c-2k}(a+k+1)_{c-2k-1}}{(a+k+\frac{1}{2})_{c-2k}(k+1)_{c-2k-1}}
\end{equation}
and
\begin{equation}
    \frac{M(H_{2a+1,c,c;d})}{M(H_{2a+1,c,c})}=\frac{1}{4^{d}}\prod_{k=0}^{d-1} \frac{(a+k+1)_{c-2k}(k+\frac{3}{2})_{c-2k-2}}{(k+1)_{c-2k-1}(a+k+\frac{3}{2})_{c-2k-1}}.  
\end{equation}

We only prove (4.1) in detail because the proof of (4.2) is completely analogous. As we saw in the previous section, by the Factorization Theorem,
\begin{equation}
    M(H_{2a,c,c;d})=2^{c-d}M(H_{2a,c,c;d}^{+})M(H_{2a,c,c;d}^{-}).
\end{equation}
One can make the following observations about the regions $H_{2a,c,c;d}^{+}$ and $H_{2a,c,c;d}^{-}$.\\
1) $M(H_{2a,c,c;d}^{-})=M(\overline{R}_{c-1,0,a})$, which is nonzero and does not depend on the parameter $d$.\\
2) $M(H_{2a,c,c;d}^{+})=
\begin{cases}
    M(R_{0,c,a-1})& \text{if } d=0,\\
    M(\overline{R}_{d-1,c-d,a})& \text{if } 1\leq d \leq c.
\end{cases}
$

After removing forced lozenges from $H_{2a,c,c;d}^{-}$, one obtains the region $\overline{R}_{c-1,0,a}$, so 1) follows. 2) can also be explained in a similar way (see Figure 4.1).

\begin{figure}
    \centering
    \includegraphics[width=10cm]{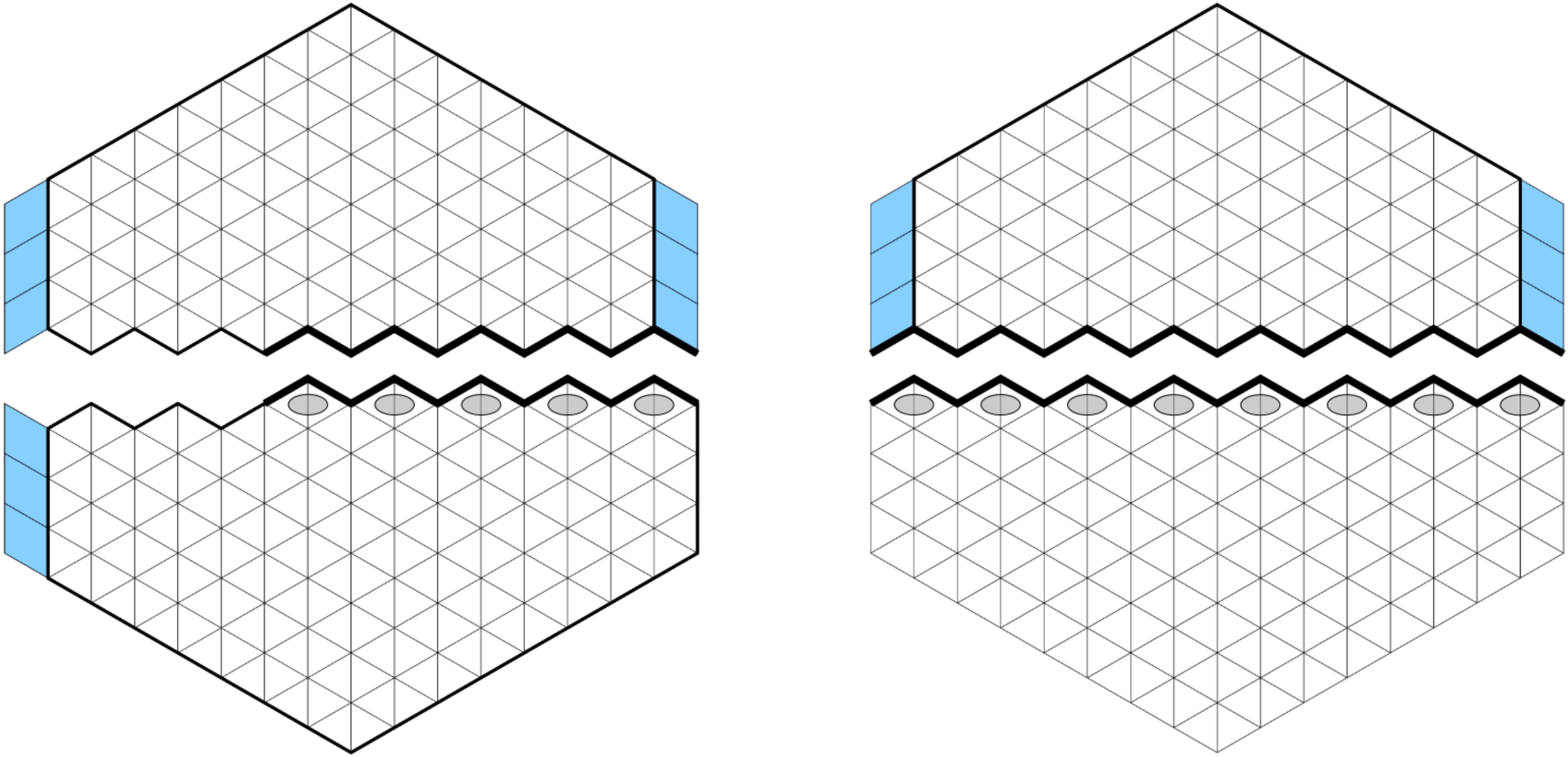}
    \caption{$H_{6,8,8;3}^+$ (bottom left), $H_{6,8,8;3}^-$ (top left), $H_{6,8,8;0}^+$ (bottom right), and $H_{6,8,8;0}^-$ (top right). Forced lozenges are indicated by shading.}
\end{figure}

Using these observations and $M(H_{2a,c,c;0})=M(H_{2a,c,c})$,

\begin{equation}
\begin{aligned}
    \frac{M(H_{2a,c,c;d})}{M(H_{2a,c,c})}=\frac{M(H_{2a,c,c;d})}{M(H_{2a,c,c;0})}=\prod_{k=0}^{d-1}\frac{M(H_{2a,c,c;k+1})}{M(H_{2a,c,c;k})}
    &=\prod_{k=0}^{d-1}\frac{2^{c-k-1}M(H_{2a,c,c;k+1}^{+})M(H_{2a,c,c;k+1}^{-})}{2^{c-k}M(H_{2a,c,c;k}^{+})M(H_{2a,c,c;k}^{-})}\\
    &=\prod_{k=0}^{d-1}\frac{M(H_{2a,c,c;k+1}^{+})}{2 M(H_{2a,c,c;k}^{+})}.
\end{aligned}
\end{equation}
Using Theorem 3.1, one can check that
\begin{equation}
    \frac{M(H_{2a,c,c;1}^{+})}{2M(H_{2a,c,c;0}^{+})}=\frac{M(\overline{R}_{0,c-1,a})}{2M(R_{0,c,a-1})}=\frac{(\frac{1}{2})_{c}(a+1)_{c-1}}{(a+\frac{1}{2})_{c}(1)_{c-1}}
\end{equation}

and

\begin{equation}
    \frac{M(H_{2a,c,c;k+1}^{+})}{2M(H_{2a,c,c;k}^{+})}=\frac{M(\overline{R}_{k,c-k-1,a})}{2M(\overline{R}_{k-1,c-k,a})}=\frac{(1)_{k}(\frac{1}{2})_{c-k}(a+k+1)_{c-2k-1}}{(1)_{c-k-1}(\frac{1}{2})_{k}(a+k+\frac{1}{2})_{c-2k}}=\frac{(k+\frac{1}{2})_{c-2k}(a+k+1)_{c-2k-1}}{(a+k+\frac{1}{2})_{c-2k}(k+1)_{c-2k-1}}
\end{equation}
for any $k$ such that $1\leq k \leq c-1$.
Thus, by combining (4.4), (4.5), and (4.6), we obtain (4.1). This completes the proof of (2.1) when $c-b=0$.

The proof of (4.2) is very similar to that of (4.1). Again, from the application of the Factorization Theorem,

\begin{equation}
    M(H_{2a+1,c,c;d})=2^{c-d}M(H_{2a+1,c,c;d}^+)M(H_{2a+1,c,c;d}^-).
\end{equation}

The following observations about $H_{2a+1,c,c;d}^{+}$ and $H_{2a+1,c,c;d}^{-}$ are needed to prove (4.2) (see Figure 4.2).\\
$1^\prime$) $M(H_{2a+1,c,c;d}^{-})=M(\overline{R}_{c,0,a})$, which is nonzero and does not depend on the parameter $d$.\\
$2^\prime$) $M(H_{2a+1,c,c;d}^{+})=
\begin{cases}
    M(R_{d,c-d-1,a})& \text{if } 0\leq d \leq c-1,\\
    M(\overline{R}_{c,0,a})& \text{if } d=c.
\end{cases}
$

Together with Theorem 3.1, one can repeat the argument presented in the proof of (4.1) and show that (2.2) holds when $c-b=0$. This completes the proof of the theorem when $c-b=0$.\\


\begin{figure}
    \centering
    \includegraphics[width=10cm]{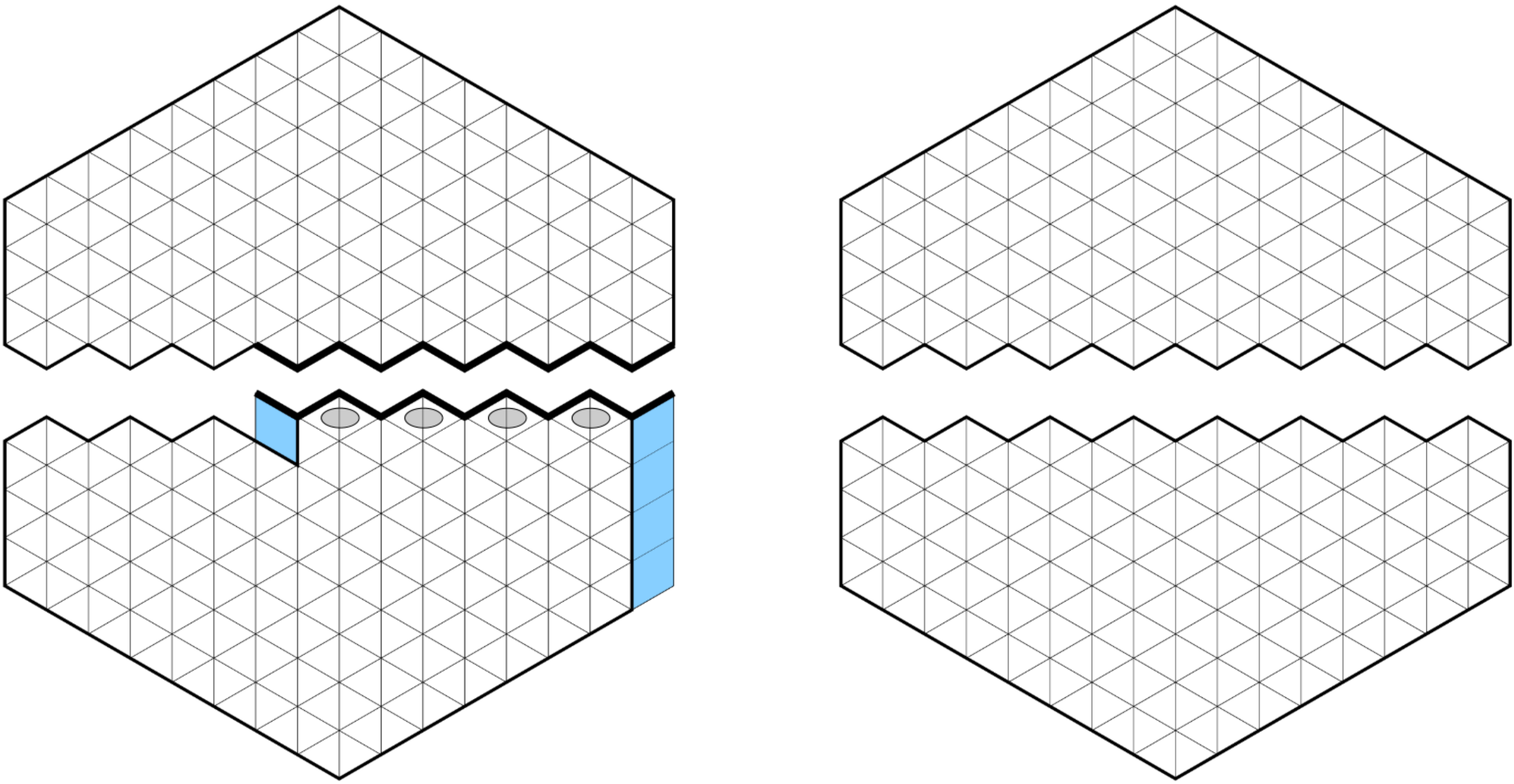}
    \caption{$H_{7,8,8;3}^+$ (bottom left), $H_{7,8,8;3}^-$ (top left), $H_{7,8,8;8}^+$ (bottom right), and $H_{7,8,8;8}^-$ (top right). Forced lozenges are indicated by shading.}
\end{figure}






\textbf{Step 2: the case } $\mathbf{c-b=1}$

Now, we make two recurrence relations involving the two regions $H_{2a,b,c;k}$ and $H_{2a+1,b,c;k}$ using the two versions of Kuo's graphical condensation mentioned in the previous section.

First, for nonnegative integers $a,k$ and positive integers $b,c$ that satisfy $k\leq \min(b,c)$, we define the two regions $H'_{2a,b,c;k}$ and $H''_{2a,b,c;k}$ as follows (see Figure 4.3 for two examples). We consider a hexagon with sides of length $2a, b+1, c, 2a+1, b$, and $c+1$ clockwise from the left. Consider the perpendicular bisector of the left side of the region, and label the unit triangles on the bisector that are also contained in the region by $1, 2, ...$ from the left. $H'_{2a,b,c;k}$ is the region obtained from the hexagonal region by deleting $2k+1$ unit triangles labeled by $1, 2, ... , 2k+1$. Similarly, $H''_{2a,b,c;k}$ is the region obtained from the same hexagonal region by deleting $2k$ unit triangles labeled by $1, 2, ... , 2k$. From the region $H'_{2a,b,c;k}$, we choose four unit triangles $x,y,z,$ and $w$ as described in the left picture of Figure 4.3. We also choose four unit triangles $x, y, z,$ and $w$ from $H''_{2a,b,c;k}$ as specified in the right picture of Figure 4.3.

\begin{figure}
    \centering
    \includegraphics[width=10cm]{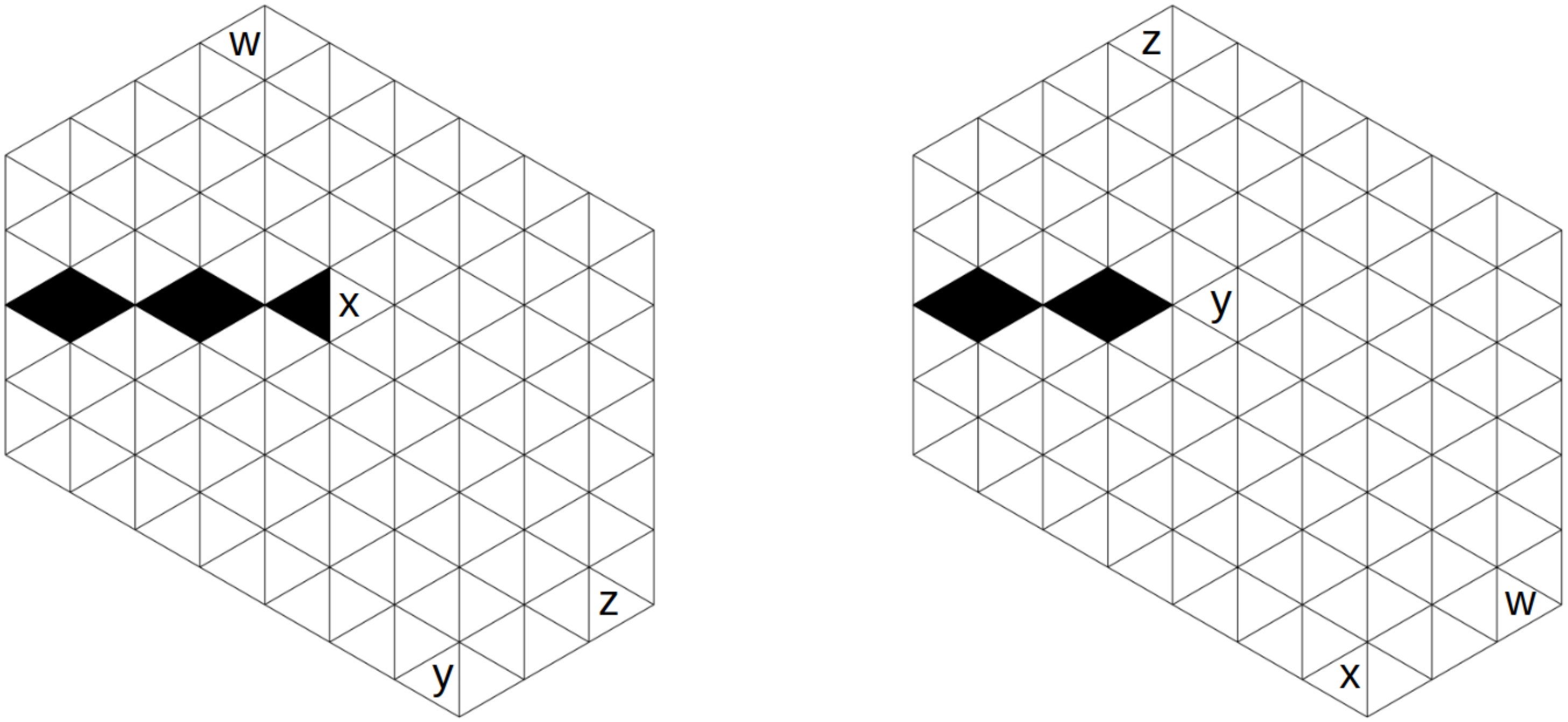}
    \caption{The regions $H'_{2a,b,c;k}$ (left) and $H''_{2a,b,c;k}$ (right) with $a=2$, $b=3$, $c=6$, and $k=2$. Positions of four unit triangles $x, y, z,$ and $w$ are specified.}
\end{figure}

One can easily see that the dual graphs of the two regions $H'_{2a,b,c;k}$ and $H''_{2a,b,c;k}$ and the choices of the four vertices satisfy the assumptions of Theorem 3.2 and Theorem 3.3, respectively. Using the first version of Kuo's graphical condensation (Theorem 3.2) on the region $H'_{2a,b,c;k}$ with four unit triangles $x,y,z,$ and $w$, we obtain the following recurrence (see Figure 4.4 that shows six regions appearing in the recurrence):

\begin{figure}
    \centering
    \includegraphics[width=13.5cm]{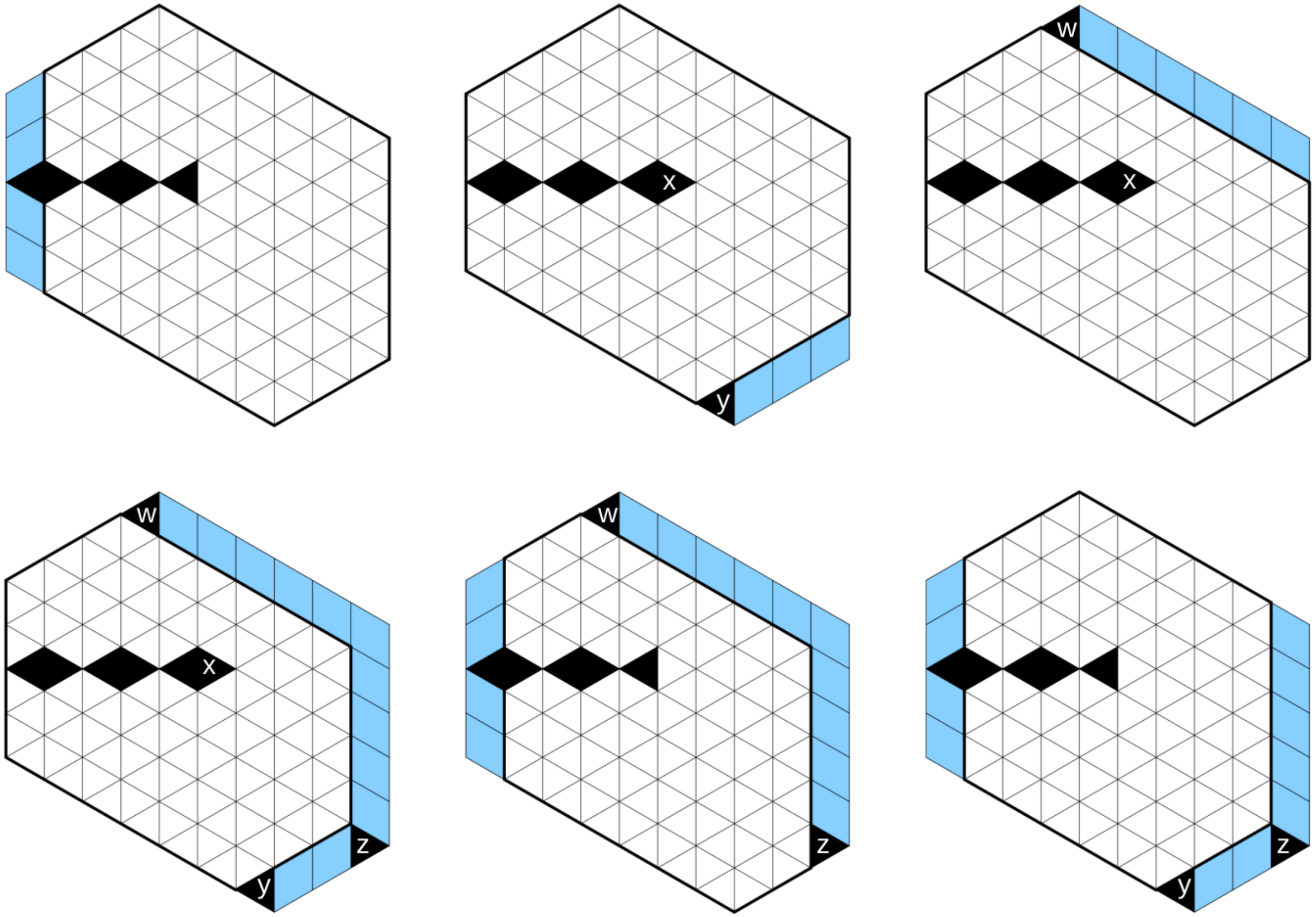}
    \caption{The six regions appearing in the recurrence relation obtained from Kuo's graphical condensation on $H'_{2a,b,c;k}$ with four unit triangles $x,y,z,$ and $w$.}
\end{figure}

\begin{equation}
\begin{aligned}
    &M(H_{2a+1,b,c;k})M(H_{2a,b,c;k+1})\\
    &=M(H_{2a,b+1,c;k+1})M(H_{2a+1,b-1,c;k})+M(H_{2a,b,c+1;k+1})M(H_{2a+1,b,c-1;k}).
\end{aligned}
\end{equation}

On the other hand, using the second version of Kuo's graphical condensation (Theorem 3.3) on the region $H''_{2a,b,c,k}$ with four unit triangles $x,y,z,$ and $w$, we obtain the following recurrence (see Figure 4.5 that shows six regions appearing in the recurrence):

\begin{figure}
    \centering
    \includegraphics[width=13.5cm]{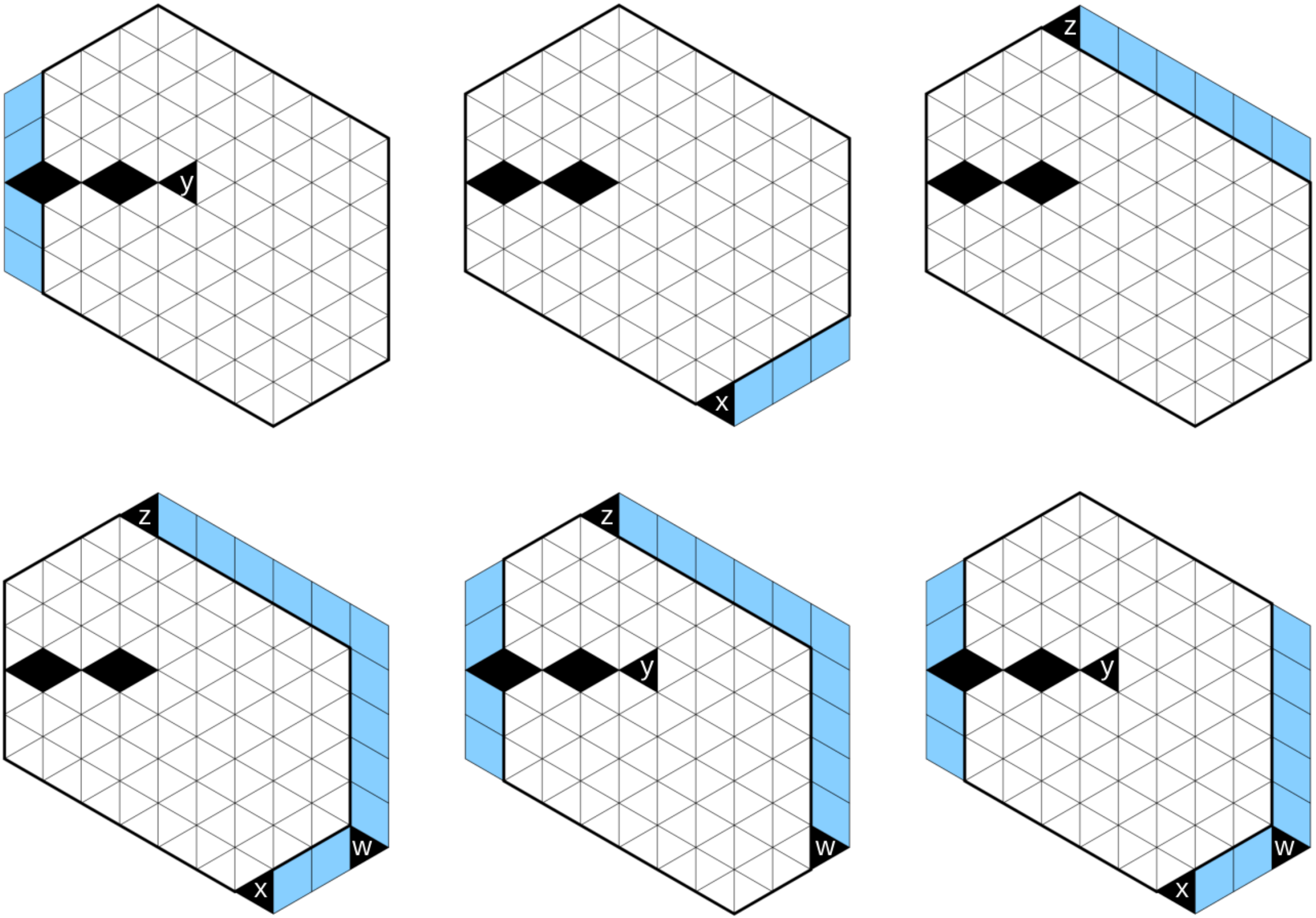}
    \caption{The six regions appearing in the recurrence relation obtained from Kuo's graphical condensation on $H''_{2a,b,c,k}$ with four unit triangles $x,y,z,$ and $w$.}
\end{figure}

\begin{equation}
\begin{aligned}
    &M(H_{2a+1,b,c;k})M(H_{2a,b,c;k})=M(H_{2a,b+1,c;k})M(H_{2a+1,b-1,c;k})+M(H_{2a,b,c+1;k})M(H_{2a+1,b,c-1;k}).
\end{aligned}
\end{equation}

We specialize (4.8) and (4.9) to prove the case when $c-b=1$.
If we replace $b$ by $c$ in (4.8), then

\begin{equation}
\begin{aligned}
    &M(H_{2a+1,c,c;k})M(H_{2a,c,c;k+1})\\
    &=M(H_{2a,c+1,c;k+1})M(H_{2a+1,c-1,c;k})+M(H_{2a,c,c+1;k+1})M(H_{2a+1,c,c-1;k}).
\end{aligned}
\end{equation}
Since the regions $H_{2a,c+1,c;k+1}$ and $H_{2a+1,c,c-1;k}$ have the same number of lozenge tilings as $H_{2a,c,c+1;k+1}$ and $H_{2a+1,c-1,c;k}$, respectively, the equation (4.10) can be rewritten as follows:
\begin{equation}
    M(H_{2a+1,c,c;k})M(H_{2a,c,c;k+1})=2M(H_{2a,c,c+1;k+1})M(H_{2a+1,c-1,c;k}).
\end{equation}

Similarly, if we replace $b$ by $c$ in (4.9), then

\begin{equation}
\begin{aligned}
    &M(H_{2a+1,c,c;k})M(H_{2a,c,c;k})=M(H_{2a,c+1,c;k})M(H_{2a+1,c-1,c;k})+M(H_{2a,c,c+1;k})M(H_{2a+1,c,c-1;k}),
\end{aligned}
\end{equation}

and the same reasoning leads us to the following equation

\begin{equation}
    M(H_{2a+1,c,c;k})M(H_{2a,c,c;k})=2M(H_{2a,c,c+1;k})M(H_{2a+1,c-1,c;k}).
\end{equation}

If we divide each side of (4.11) by that of (4.13)\footnote{We can do this division because the left side of (4.13) is nonzero because of the case when $c-b=0$.}, then we obtain

\begin{equation}
    \frac{M(H_{2a,c,c;k+1})}{M(H_{2a,c,c;k})}=\frac{M(H_{2a,c,c+1;k+1})}{M(H_{2a,c,c+1;k})}.
\end{equation}

Thus, by (4.1) and (4.14), for any $d$ such that $0\leq d \leq c$,

\begin{equation}
\begin{aligned}
    \frac{M(H_{2a,c,c+1;d})}{M(H_{2a,c,c+1})}=\frac{M(H_{2a,c,c+1;d})}{M(H_{2a,c,c+1;0})}=\prod_{k=0}^{d-1}\frac{M(H_{2a,c,c+1;k+1})}{M(H_{2a,c,c+1;k})}&=\prod_{k=0}^{d-1}\frac{M(H_{2a,c,c;k+1})}{M(H_{2a,c,c;k})}\\
    &=\prod_{k=0}^{d-1} \frac{(k+\frac{1}{2})_{c-2k}(a+k+1)_{c-2k-1}}{(k+1)_{c-2k-1}(a+k+\frac{1}{2})_{c-2k}}.
\end{aligned}
\end{equation}
This proves (2.1) when $c-b=1$. To prove (2.2), we replace $k$ by $k+1$ in (4.13) and divide each side of the obtained equation by that of (4.11)\footnote{We can do this division because the left side of (4.11) is nonzero because of the case when $c-b=0$.}, thereby obtaining

\begin{equation}
    \frac{M(H_{2a+1,c,c;k+1})}{M(H_{2a+1,c,c;k})}=\frac{M(H_{2a+1,c-1,c;k+1})}{M(H_{2a+1,c-1,c;k})}.
\end{equation}

Hence, by (4.2) and (4.16), for any $d$ such that $0\leq d \leq c-1$,
\begin{equation}
\begin{aligned}
    \frac{M(H_{2a+1,c-1,c;d})}{M(H_{2a+1,c-1,c})}=\frac{M(H_{2a+1,c-1,c;d})}{M(H_{2a+1,c-1,c;0})}&=\prod_{k=0}^{d-1}\frac{M(H_{2a+1,c-1,c;k+1})}{M(H_{2a+1,c-1,c;k})}\\
    &=\prod_{k=0}^{d-1}\frac{M(H_{2a+1,c,c;k+1})}{M(H_{2a+1,c,c;k})} \\
    &=\frac{1}{4^{d}}\prod_{k=0}^{d-1} \frac{(a+k+1)_{c-2k}(k+\frac{3}{2})_{c-2k-2}}{(k+1)_{c-2k-1}(a+k+\frac{3}{2})_{c-2k-1}}.
\end{aligned}
\end{equation}
This proves that (2.2) holds when $c-b=1$ and completes the proof of the theorem when $c-b=1$.\\

\textbf{Step 3: the case } $\mathbf{c-b\geq 2}$

Now, we show (2.1) and (2.2) hold in general. We prove them using induction on the value $c-b$ (we call it \textit{outer induction}). We already verified that (2.1) and (2.2) are true when $c-b=0$ or $1$. Suppose that (2.1) and (2.2) hold when $c-b<s$ for some $s \geq 2$. Under this assumption, we need to verify that (2.1) and (2.2) still hold when $c-b=s$.

We prove this induction step using another induction on $d$ (we call it \textit{inner induction}). Since the regions $H_{2a,b,c;0}$ and $H_{2a+1,b,c;0}$ are the same as $H_{2a,b,c}$ and $H_{2a+1,b,c}$, the theorem holds when $c-b=s$ and $d=0$. Suppose that (2.1) and (2.2) are true when $c-b=s$ and $d < t$ for some $1 \leq t \leq b$. Under these assumptions, we have to show that (2.1) and (2.2) still hold when $c-b=s$ and $d=t$. If we replace $c$ by $c-1$ and $k$ by $t-1$ in (4.8), then we obtain

\begin{equation}
\begin{aligned}
    &M(H_{2a+1,b,c-1;t-1})M(H_{2a,b,c-1;t})\\
    &=M(H_{2a,b+1,c-1;t})M(H_{2a+1,b-1,c-1;t-1})+M(H_{2a,b,c;t})M(H_{2a+1,b,c-2;t-1}).
\end{aligned}
\end{equation}

In (4.18), the two terms $M(H_{2a,b,c-1;t})$ and $M(H_{2a,b+1,c-1;t})$ are given by (2.1) by the induction hypothesis of the outer induction. Also, the other two terms $M(H_{2a+1,b,c-1;t-1})$ and $M(H_{2a+1,b,c-2;t-1})$ are given by (2.2), again by the induction hypothesis of the outer induction.
Furthermore, the term $M(H_{2a+1,b-1,c-1;t-1})$ is given by (2.2) by the induction hypothesis of the inner induction.

Hence, to show that $M(H_{2a,b,c;t})$ is given by (2.1), it is enough to show that the following identity

\begin{equation}
    \frac{M(H_{2a,b+1,c-1;t})M(H_{2a+1,b-1,c-1;t-1})}{M(H_{2a+1,b,c-1;t-1})M(H_{2a,b,c-1;t})}+\frac{M(H_{2a,b,c;t})M(H_{2a+1,b,c-2;t-1})}{M(H_{2a+1,b,c-1;t-1})M(H_{2a,b,c-1;t})}=1
\end{equation}
is true when we replace every $M(\cdot)$ in (4.19) by the formulas in (2.1) and (2.2) together with MacMahon's formula.

Using (2.1), (2.2), and MacMahon's formula, one can easily check that the following identities hold:

\begin{equation}
    \frac{M(H_{2a,b+1,c-1;t})}{M(H_{2a,b,c-1;t})}=\frac{b!(2a+b+c-1)!}{(2a+b)!(b+c-1)!}\prod_{k=0}^{t-1}\Bigg[\frac{(a+b-k)(b+c-2k-1)}{(b-k)(2a+b+c-2k-1)}\Bigg],
\end{equation}

\begin{equation}
    \frac{M(H_{2a+1,b-1,c-1;t-1})}{M(H_{2a+1,b,c-1;t-1})}=\frac{(2a+b)!(b+c-2)!}{(2a+b+c-1)!(b-1)!}\prod_{k=0}^{t-2}\Bigg[\frac{(b-k-1)(2a+b+c-2k-1)}{(a+b-k)(b+c-2k-3)}\Bigg],
\end{equation}

\begin{equation}
    \frac{M(H_{2a,b,c;t})}{M(H_{2a,b,c-1;t})}=\frac{(c-1)!(2a+b+c-1)!}{(2a+c-1)!(b+c-1)!}\prod_{k=0}^{t-1}\Bigg[\frac{(a+c-k-1)(b+c-2k-1)}{(c-k-1)(2a+b+c-2k-1)}\Bigg],
\end{equation}
and
\begin{equation}
    \frac{M(H_{2a+1,b,c-2;t-1})}{M(H_{2a+1,b,c-1;t-1})}=\frac{(2a+c-1)!(b+c-2)!}{(2a+b+c-1)!(c-2)!}\prod_{k=0}^{t-2}\Bigg[\frac{(c-k-2)(2a+b+c-2k-1)}{(a+c-k-1)(b+c-2k-3)}\Bigg].
\end{equation}

If we put (4.20)-(4.23) in the left side of (4.19), then

\begin{equation}
\begin{aligned}
    &\frac{M(H_{2a,b+1,c-1;t})M(H_{2a+1,b-1,c-1;t-1})}{M(H_{2a+1,b,c-1;t-1})M(H_{2a,b,c-1;t})}+\frac{M(H_{2a,b,c;t})M(H_{2a+1,b,c-2;t-1})}{M(H_{2a+1,b,c-1;t-1})M(H_{2a,b,c-1;t})}\\
    &=\frac{M(H_{2a,b+1,c-1;t})}{M(H_{2a,b,c-1;t})}\cdot\frac{M(H_{2a+1,b-1,c-1;t-1})}{M(H_{2a+1,b,c-1;t-1})}+\frac{M(H_{2a,b,c;t})}{M(H_{2a,b,c-1;t})}\cdot\frac{M(H_{2a+1,b,c-2;t-1})}{M(H_{2a+1,b,c-1;t-1})}\\
    &=\frac{b}{b+c-1}\cdot\frac{(a+b-t+1)(b+c-1)}{b(2a+b+c-2t+1)}+\frac{c-1}{b+c-1}\cdot\frac{(a+c-t)(b+c-1)}{(c-1)(2a+b+c-2t+1)}\\
    &=\frac{a+b-t+1}{2a+b+c-2t+1}+\frac{a+c-t}{2a+b+c-2t+1}\\
    &=1.
\end{aligned}
\end{equation}

Hence, we have shown that $M(H_{2a,b,c;t})$ is given by (2.1). We continue to show that $M(H_{2a+1,b,c;t})$ is also given by (2.2) under the same induction hypothesis. If we replace $b$ by $b+1$ and $k$ by $t$ in (4.9), then we obtain

\begin{equation}
\begin{aligned}
    &M(H_{2a+1,b+1,c;t})M(H_{2a,b+1,c;t})\\
    &=M(H_{2a,b+2,c;t})M(H_{2a+1,b,c;t})+M(H_{2a,b+1,c+1;t})M(H_{2a+1,b+1,c-1;t}).
\end{aligned}
\end{equation}

In (4.25), the two terms $M(H_{2a,b+1,c;t})$ and $M(H_{2a,b+2,c;t})$ are given by (2.1) by the induction hypothesis of the outer induction. Also, the other two terms $M(H_{2a+1,b+1,c;t})$ and $M(H_{2a+1,b+1,c-1;t})$ are given by (2.2), again by the induction hypothesis of the outer induction.
Also, we already showed that the term $M(H_{2a,b+1,c+1;t})$ is given by (2.1) under the same hypothesis.

Hence, to show that $M(H_{2a+1,b,c;t})$ is given by (2.2), it is enough to show that the following identity

\begin{equation}
    \frac{M(H_{2a,b+2,c;t})M(H_{2a+1,b,c;t})}{M(H_{2a+1,b+1,c;t})M(H_{2a,b+1,c;t})}+\frac{M(H_{2a,b+1,c+1;t})M(H_{2a+1,b+1,c-1;t})}{M(H_{2a+1,b+1,c;t})M(H_{2a,b+1,c;t})}=1
\end{equation}
holds when we replace every $M(\cdot)$ by the formulas in (2.1) and (2.2) together with MacMahon's formula.

Again, using (2.1), (2.2), and MacMahon's formula, it is straightforward to check that the following identities are true:

\begin{equation}
    \frac{M(H_{2a,b+2,c;t})}{M(H_{2a,b+1,c;t})}=\frac{(b+1)!(2a+b+c+1)!}{(b+c+1)!(2a+b+1)!}\prod_{k=0}^{t-1}\frac{(a+b-k+1)(b+c-2k+1)}{(b-k+1)(2a+b+c-2k+1)},
\end{equation}

\begin{equation}
    \frac{M(H_{2a+1,b,c;t})}{M(H_{2a+1,b+1,c;t})}=\frac{(b+c)!(2a+b+1)!}{b!(2a+b+c+1)!}\prod_{k=0}^{t-1}\frac{(b-k)(2a+b+c-2k+1)}{(a+b-k+1)(b+c-2k-1)},
\end{equation}

\begin{equation}
    \frac{M(H_{2a,b+1,c+1;t})}{M(H_{2a,b+1,c;t})}=\frac{c!(2a+b+c+1)!}{(b+c+1)!(2a+c)!}\prod_{k=0}^{t-1}\frac{(a+c-k)(b+c-2k+1)}{(c-k)(2a+b+c-2k+1)},
\end{equation}
and
\begin{equation}
    \frac{M(H_{2a+1,b+1,c-1;t})}{M(H_{2a+1,b+1,c;t})}=\frac{(b+c)!(2a+c)!}{(c-1)!(2a+b+c+1)!}\prod_{k=0}^{t-1}\frac{(c-k-1)(2a+b+c-2k+1)}{(a+c-k)(b+c-2k-1)}.
\end{equation}

If we put (4.27)-(4.30) in the left side of (4.26), then

\begin{equation}
\begin{aligned}
    &\frac{M(H_{2a,b+2,c;t})M(H_{2a+1,b,c;t})}{M(H_{2a+1,b+1,c;t})M(H_{2a,b+1,c;t})}+\frac{M(H_{2a,b+1,c+1;t})M(H_{2a+1,b+1,c-1;t})}{M(H_{2a+1,b+1,c;t})M(H_{2a,b+1,c;t})}\\
    &\frac{M(H_{2a,b+2,c;t})}{M(H_{2a,b+1,c;t})}\cdot\frac{M(H_{2a+1,b,c;t})}{M(H_{2a+1,b+1,c;t})}+\frac{M(H_{2a,b+1,c+1;t})}{M(H_{2a,b+1,c;t})}\cdot\frac{M(H_{2a+1,b+1,c-1;t})}{M(H_{2a+1,b+1,c;t})}\\
    &=\frac{b+1}{b+c+1}\cdot\frac{(b-t+1)(b+c+1)}{(b+1)(b+c-2t+1)}+\frac{c}{b+c+1}\cdot\frac{(c-t)(b+c+1)}{c(b+c-2t+1)}\\
    &=\frac{b-t+1}{b+c-2t+1}+\frac{c-t}{b+c-2t+1}\\
    &=1.
\end{aligned}
\end{equation}

Hence, we have shown that $M(H_{2a+1,b,c;t})$ is given by (2.2), and this verifies the induction step for the inner induction.
Thus, (2.1) and (2.2) are true for any $b,c,$ and $d$ such that $c-b=s$ and $0\leq d \leq b$, and this implies that the induction step for the outer induction is also verified.
Hence, by mathematical induction, we can conclude that (2.1) and (2.2) are true for any $b$ and $c$ such that $c-b\geq0$.
This completes the proof of the theorem.
\begin{flushright}
\qedsymbol
\end{flushright}

\section{A proof of Theorem 2.3}

Since the proof of (2.6) and that of (2.7) are very similar, we only present the proof of (2.6).

By Theorem 2.1,

\begin{equation}
\begin{aligned}
    &\frac{M(H_{2aN,bN,cN;dN})}{M(H_{2aN,bN,cN})}\\
    &=\prod_{k=0}^{dN-1} \frac{(k+\frac{1}{2})_{bN-2k}(aN+k+1)_{bN-2k-1}(bN-k+\frac{1}{2})_{\frac{cN-bN}{2}}(cN-k)_{-\frac{cN-bN}{2}}}{(aN+k+\frac{1}{2})_{bN-2k}(k+1)_{bN-2k-1}(aN+bN-k+\frac{1}{2})_{\frac{cN-bN}{2}}(aN+cN-k)_{-\frac{cN-bN}{2}}}.
\end{aligned}
\end{equation}

We use the following identities, which can be easily verified from the (extended) definition of the shifted factorial. For positive integer $i$ and an integer $j$ such that $i+j>0$ holds,
\begin{equation}
    (i)_{j}=\frac{(i+j-1)!}{(i-1)!},
\end{equation}
and
\begin{equation}
    \Big(i+\frac{1}{2}\Big)_{j}=\frac{1}{2^{2j}}\frac{i!(2i+2j)!}{(i+j)!(2i)!}=\frac{1}{2^{2j}}\frac{(i-1)!(2i+2j-1)!}{(i+j-1)!(2i-1)!}.
\end{equation}

By applying (5.2) and (5.3) to (5.1), we obtain

\begin{equation}
\begin{aligned}
    \frac{M(H_{2aN,bN,cN,dN})}{M(H_{2aN,bN,cN})}
    =\prod_{k=0}^{dN-1} &\Bigg[\frac{k!}{(aN+k)!}\frac{k!}{(aN+k)!}\frac{(2aN+2k)!}{(2k)!}\frac{(aN+bN-k-1)!}{(bN-k-1)!}\\
    &\times\frac{(aN+cN-k-1)!}{(cN-k-1)!}\frac{(bN+cN-2k-1)!}{(2aN+bN+cN-2k-1)!}\Bigg].
\end{aligned}
\end{equation}

Since a product of factorials of consecutive integers can be expressed as a ratio of the Barnes G-functions,
\begin{equation}
    \prod_{k=0}^{dN-1} \frac{k!}{(aN+k)!} = \frac{G(dN+1)G(aN+1)}{G((a+d)N+1)},
\end{equation}
\begin{equation}
    \prod_{k=0}^{dN-1} \frac{(aN+bN-k-1)!}{(bN-k-1)!} = \frac{G((a+b)N+1)G((b-d)N+1)}{G((a+b-d)N+1)G(bN+1)},
\end{equation}
and
\begin{equation}
    \prod_{k=0}^{dN-1} \frac{(aN+cN-k-1)!}{(cN-k-1)!} = \frac{G((a+c)N+1)G((c-d)N+1)}{G((a+c-d)N+1)G(cN+1)}.
\end{equation}

The remaining two factors
\begin{equation*}
    \displaystyle \prod_{k=0}^{dN-1} \frac{(2aN+2k)!}{(2k)!}    
\end{equation*}
and
\begin{equation*}
    \displaystyle \prod_{k=0}^{dN-1} \frac{(bN+cN-2k-1)!}{(2aN+bN+cN-2k-1)!}    
\end{equation*}
can be also expressed in terms of the Barnes G-functions and factorials. For example,
\begin{equation}
\begin{aligned}
    \prod_{k=0}^{dN-1} \frac{(2aN+2k)!}{(2k)!}&=\Bigg(\Bigg[\prod_{k=0}^{dN-1} \frac{(2aN+2k)!}{(2k)!}\Bigg]\Bigg[\prod_{k=0}^{dN-1} \frac{(2aN+2k+1)!}{(2k+1)!}\Bigg]\Bigg[\prod_{k=0}^{dN-1} \frac{2k+1}{2aN+2k+1}\Bigg]\Bigg)^{\frac{1}{2}}\\
    &=\Bigg(\Bigg[\prod_{k=0}^{2dN-1} \frac{(2aN+k)!}{(k)!}\Bigg]\Bigg[\prod_{k=0}^{dN-1} \frac{2k+1}{2aN+2k+1}\Bigg]\Bigg)^{\frac{1}{2}}\\
    &=\Bigg[\frac{G((2a+2d)N+1)}{G(2aN+1)G(2dN+1)}\frac{(2dN)!((a+d)N)!(2aN)!}{(dN)!(aN)!((2a+2d)N)!}\Bigg]^{\frac{1}{2}}
\end{aligned}
\end{equation}
and similarly,
\begin{equation}
\begin{aligned}
    &\prod_{k=0}^{dN-1} \frac{(bN+cN-2k-1)!}{(2aN+bN+cN-2k-1)!}\\
    =&\Bigg[\frac{G((b+c)N+1)G((2a+b+c-2d)N+1)}{G((2a+b+c)N+1)G((b+c-2d)N+1)}\\
    &\times\frac{((b+c)N)!((2a+b+c-2d)N)!((a+\frac{1}{2}b+\frac{1}{2}c)N)!((\frac{1}{2}b+\frac{1}{2}c-d)N)!}{((b+c-2d)N)!((2a+b+c)N)!((a+\frac{1}{2}b+\frac{1}{2}c-d)N)!((\frac{1}{2}b+\frac{1}{2}c)N)!}\Bigg]^{\frac{1}{2}}.
\end{aligned}
\end{equation}

By putting (5.5)-(5.9) in (5.4), we are able to express $\displaystyle \frac{M(H_{2aN,bN,cN;dN})}{M(H_{2aN,bN,cN})}$ using only the Barnes G-functions and factorials.
Finally, using Stirling's formula and the definition of the Glaisher-Kinkelin constant\footnote{$G(n+1) \sim A^{-1}(2\pi)^{\frac{n}{2}}n^{\frac{1}{2}n^2-\frac{1}{12}}e^{-\frac{3}{4}n^{2}+\frac{1}{12}}$ as $n\to\infty$. Here, $A$ is the Glaisher-Kinkelin constant.}, we obtain the expression given on the right side of (2.6). This completes the proof.\begin{flushright}
\qedsymbol
\end{flushright}

\section{Concluding Remarks}

In this paper, we provide the product formulas for the number of lozenge tilings of two regions $H_{2a,b,c;d}$ and $H_{2a+1,b,c;d}$. One can easily see that our formulas (2.1) and (2.2) are equivalent to the following equations involving ratios of the number of lozenge tilings of two closely related regions:
\begin{equation}
    \frac{M(H_{2a,b,c;k+1})}{M(H_{2a,b,c;k})}=\frac{(k+\frac{1}{2})_{b-2k}(a+k+1)_{b-2k-1}(b-k+\frac{1}{2})_{\lfloor\frac{c-b}{2}\rfloor}(c-k)_{-\lfloor\frac{c-b}{2}\rfloor}}{(a+k+\frac{1}{2})_{b-2k}(k+1)_{b-2k-1}(a+b-k+\frac{1}{2})_{\lfloor\frac{c-b}{2}\rfloor}(a+c-k)_{-\lfloor\frac{c-b}{2}\rfloor}}
\end{equation}
and
\begin{equation}
    \frac{M(H_{2a+1,b,c;k+1})}{M(H_{2a+1,b,c;k})}=\frac{1}{4}\frac{(a+k+1)_{c-2k}(k+\frac{3}{2})_{c-2k-2}(b-k)_{\lfloor\frac{c-b}{2}\rfloor}(c-k-\frac{1}{2})_{-\lfloor\frac{c-b}{2}\rfloor}}{(k+1)_{c-2k-1}(a+k+\frac{3}{2})_{c-2k-1}(a+b-k+1)_{\lfloor\frac{c-b}{2}\rfloor}(a+c-k+\frac{1}{2})_{-\lfloor\frac{c-b}{2}\rfloor}}.
\end{equation}


The simplicity of (6.1) and (6.2) calls for a more simple proof. It would be interesting to find a direct proof of the equations (6.1) and (6.2), which will lead us to a new proof of our main theorem.\\

\textbf{Acknowledgments.}
The author thanks his advisor Professor Mihai Ciucu for his continuing encouragement and motivation. The author also thanks anonymous reviewers for carefully reading the original version of the paper and giving helpful comments. David Wilson's program vaxmacs was helpful when the author tried to find the formulas.

\end{document}